\documentclass[12pt,reqno,
]{amsart}

\usepackage[arrow,matrix,curve]{xy}

\usepackage[dvips]{graphicx} 

\usepackage{amssymb, latexsym, amsmath, amscd, array,
}

\newcommand{\hid} [1] {}

\hid{
\makeatletter
\if@twoside%
\addtolength{\evensidemargin}{20.0ex} 
\addtolength{\oddsidemargin}{-23.0ex} 
\else\fi%
\makeatother
}

\newcommand{\dR} {\mathord{^\ast\hspace*{-0.1ex}\R}}
\newcommand{\bL} {\mathbf L}

\newtheorem{theorem}{Theorem}[section]

\newtheorem{corollary}[theorem]{Corollary}

\newtheorem{construction}[theorem]{Construction}

\theoremstyle{definition}

\newtheorem{remark}[theorem]{Remark}

\numberwithin{equation}{section}

\newcommand\N {{\mathbb N}} 

\newcommand\R {{\mathbb R}} 

\newcommand\IR {{}^{*}{\mathbb R}}

\newcommand\C {{\mathbb C}}

\newcommand\Q {{\mathbb Q}} \newcommand\Z
{{\mathbb Z}}

\newcommand\Los{{\L}o{\'s}}

\newcommand\Schutzenberger{Sch\"utzen\-ber\-ger}

\newcommand{\som} {\mathcal{S}}
\newcommand{\somp} {\som'}

\newcommand\parisotes{{$\pi\alpha\rho\iota\sigma\acute{o}\tau\eta\varsigma$}}

\begin{document}

\thispagestyle{empty}

\title[Tools, objects, and chimeras] {Tools, objects, and chimeras:
Connes on the role of hyperreals in mathematics}

\author{Vladimir Kanovei} 
\address{IPPI, Moscow, and MIIT, Moscow}
\email{kanovei@rambler.ru}

\author{Mikhail G. Katz} \address{Department of Mathematics, Bar Ilan
University, Ramat Gan 52900 Israel}\email{katzmik@macs.biu.ac.il}

\author{Thomas Mormann} \address{Department of Logic and Philosophy of
Science, University of the Basque Country UPV/EHU, 20080 Donostia San
Sebastian, Spain} \email{ylxmomot@sf.ehu.es}



\begin{abstract}
We examine some of Connes' criticisms of Robinson's infinitesimals
starting in 1995.  Connes sought to exploit the Solovay
model~$\mathcal{S}$ as ammunition against non-standard analysis, but
the model tends to boomerang, undercutting Connes' own earlier work in
functional analysis.  Connes described the hyperreals as both a
``virtual theory" and a ``chimera'', yet acknowledged that his
argument relies on the transfer principle.  We analyze Connes'
``dart-throwing'' thought experiment, but reach an opposite
conclusion.  In~$\mathcal{S}$, all definable sets of reals are
Lebesgue measurable, suggesting that Connes views a theory as being
``virtual'' if it is not \emph{definable} in a suitable model of ZFC.
If so, Connes' claim that a theory of the hyperreals is ``virtual" is
refuted by the existence of a definable model of the hyperreal field
due to Kanovei and Shelah.  Free ultrafilters aren't definable, yet
Connes exploited such ultrafilters both in his own earlier work on the
classification of factors in the 1970s and 80s, and in
\emph{Noncommutative Geometry}, raising the question whether the
latter may not be vulnerable to Connes' criticism of virtuality.  We
analyze the philosophical underpinnings of Connes' argument based on
G\"odel's incompleteness theorem, and detect an apparent circularity
in Connes' logic.  We document the reliance on non-constructive
foundational material, and specifically on the Dixmier
trace~$-\hskip-9pt\int$ (featured on the front cover of Connes'
\emph{magnum opus}) and the Hahn--Banach theorem, in Connes' own
framework.  We also note an inaccuracy in Machover's critique of
infinitesimal-based pedagogy.
\end{abstract}

\maketitle

{\small
\tableofcontents
}

\section{Infinitesimals from Robinson to Connes via Choquet}

A
%
%
theory of infinitesimals claiming to vindicate Leibniz's calculus was
developed by Abraham Robinson in the 1960s (see \cite{Ro66}).  In
France, Robinson's lead was followed by G.~Reeb, G.~Choquet,%
\footnote{See e.g., Choquet's work on ultrafilters \cite{Ch}.
Choquet's constructions were employed and extended by
G.~Mokobodzki~\cite{Mo}.}
and others.  Alain Connes started his work under Choquet's leadership,
and published two texts on the hyperreals and ultrapowers (Connes
\cite{Co70a}, \cite{Co70b}).

In 1976, Connes used ultraproducts (exploiting in particular free
ultrafilters on~$\N$) in an essential manner in his work on the
classification of factors (Connes 1976, \cite{Co76}). (See
Remark~\ref{von} for Connes' use of ultrafilters in
\emph{Noncommutative geometry}.)

During the 1970s, Connes reportedly discovered that Robinson's
infinitesimals were not suitable for Connes' framework.  A quarter of
a century later, in 1995, Connes unveiled an alternative theory of
infinitesimals (Connes \cite{Co95}).  Connes' presentation of his
theory is usually not accompanied by acknowledgment of an intellectual
debt to Robinson.  Instead, it is frequently accompanied by criticism
of Robinson's framework, exploiting epithets that range from
``inadequate'' to ``end of the rope for being `explicit'\,'' (see
Table~\ref{epithets} in Section~\ref{two}).  We will examine some of
Connes' criticisms, which tend to be at tension with Connes' earlier
work.  A related challenge to the hyperreal approach was analyzed by
F.~Herzberg~\cite{He}.  Another challenge by E.~Bishop was analyzed by
Katz \& Katz \cite{KK11d}, \cite{KK11a}.  For a related analysis see
Katz \& Leichtnam \cite{KL}.

In Section~\ref{otte}, we examine the philosophical underpinnings of
Connes' position.  In Section~\ref{two}, we analyze the Connes
character and its relation to ultrafilters, and present a chronology
of Connes' criticisms of NSA.  In Section~\ref{definable}, we examine
some meta-mathematical implications of the definable model of the
hyperreal field constructed by Kanovei and Shelah.  Machover's
critique is analyzed in Section~\ref{machover}.  The power of the
\Los-Robinson transfer principle is sized up in Section~\ref{seven}.
The foundational status of the Dixmier trace and its role in
noncommutative geometry are analyzed in Section~\ref{dix}.

\section{Tools and objects}
\label{otte}

Connes' variety of Platonism can be characterized more specifically as
a \emph{prescriptive} Platonism, whereby one not merely postulates the
existence of abstract objects, but proceeds to assign ``hierarchical
levels" (see Connes \cite[p.~31]{CLS}) of realness to them, and to
issue value judgments based on the latter.  Thus, non-standard numbers
and Jordan algebras get flunking scores (see Section~\ref{shift}).
Connes mentions such ``hierarchical levels" in the context of a
dichotomy between ``tool" and ``object".  In Connes' view, only
\emph{objects} enjoy a full Platonic existence, while \emph{tools}
(such as ultrafilters and non-standard numbers) serve merely the
purpose of investigating the properties of the objects.

As a general methodological comment, we note the following.  There is
indisputably a kind of aprioriness about the natural numbers and other
concepts in mathematics, that is not accounted for by a ``formalist''
view of mathematics as a game of pushing symbols around.  Such
aprioriness requires explanation.  However, Platonism and Formalism
are not the only games in town, which is a point we will return to at
the end of the section.

To take a historical perspective on this issue, Leibniz sometimes
described infinitesimals as ``useful fictions", similar to imaginary
numbers (see Katz \& Sherry \cite{KS1}, \cite{KS2} for more details).
Leibniz's take on infinitesimals was a big novelty at the time and in
fact displeased his disciples Bernoulli, l'H\^opital, and Varignon.
But Leibniz, while clearly rejecting what would be later called a
platonist view, certainly did not think of mathematics as a
meaningless game of symbols.  One can criticize certain forms of
Platonism while adhering to the proposition that mathematics has
meaning.

\subsection{Tool/object dichotomy}

Connes' approach to the tool/object dichotomy is problematic, first
and foremost, because it does not do justice to the real history of
mathematics.  Mathematical concepts may start their career as mere
tools or instruments for manipulating concepts already given or
accepted as full-fledged objects, but later they (the tools) may
themselves become recognized as full-fledged objects.  Historical
examples of such processes abound.  The ancient Greeks did not think
of the rationals as numbers, but rather as relations among natural
numbers (see e.g., B\l aszczyk et al.~\cite[Section~2.1]{BKS}).
Wallis and others in the 17th century were struggling with the
ontological expansion involved in incorporating irrational
(transcendental) numbers beyond the algebraic ones in the number
system.  Ideal points and ideal lines at infinity in projective
geometry had to face an uphill battle before joining the ranks of
\emph{objects} that can be mentioned in ontologically polite company
(see e.g., M.~Wilson \cite{Wil92}).  G.~Cantor's cardinals started as
indices and notational subscripts for sets, and only gradually came to
be thought of as objects in their own right. Certain well-established
objects still bear the name \emph{imaginary} because they were once
characterized as not possessing the same reality as genuine objects.
R.~Hersh \cite[p.~74]{He97} describes some striking cases, including
Fourier analysis, of a historical evolution of tools into objects.

The distinction between ``tools'' and ``real objects'' is not only
blurred by the ongoing conceptual evolution of mathematics.  It is
also \emph{relative} to the perspective one takes.  For instance,
set-theoretic topology considers points as the basic building blocks
of its objects, to wit, topological spaces.  From this perspective,
nothing is a more robust and solid object than a point.  On the other
hand, from the perspective of ``point-free" (lattice-theoretical)
topology, the points of set-theoretic topology appear as highly
``chimerical" entities the existence of which can only be ensured by
relying on the axiom of choice or some similar lofty principle
(cf. Gierz et al.~\cite{Gi}).  More precisely, the situation can be
described as follows.  The basic objects of point-free topology are
complete Heyting algebras (locales) which correspond to the Heyting
algebras of open sets of topological spaces.  The prime elements of
these algebras may be considered as their ``points".  The existence of
sufficiently many points can only be secured by relying on the
Hausdorff maximality principle.  Under some mild assumptions on the
Heyting algebras and the topological spaces involved, one can show
that there is a 1-1 correspondence between set-theoretical points of
spaces and constructed points of the corresponding Heyting algebras
(cf.~ibid., Proposition V-5.20, p.~423).

\subsection{The results of Solovay and Shelah}

The perspectival relativity of the tool/object distinction and the
mutual dependence between its components do not pose a problem for an
account that recognizes both tools and objects as \emph{complementary}
components of mathematics (that would perhaps make both of them
``primordial'' in Connes' terminology; see Subsection~\ref{103}).

This may be elaborated as follows.  As in any other realm of
knowledge, also in mathematics, object and tool of knowledge are
connected through the activity of mathematical research and
application: the one does not make sense without the other.  The
dynamics of knowledge requires that both components are not only
related, but also opposed to each other.  Objects are, as the
etymological roots of this word reveal, ``resistances" or ``obstacles"
for knowledge (similarly for the Greek \emph{problema} and the German
\emph{Gegenstand}).  Tools should therefore not be disparaged as mere
subjective ``chimeras" but should be conceived of, together with
objects, as constitutive ingredients of the evolution of mathematical
knowledge (cf.~Otte 1994 \cite[ch.~X]{Ot}).%
\footnote{In a related vein, J.-P.~Marquis \cite{Ma97}, \cite{Ma06}
pointed out the ever-growing importance of complex conceptual tools
for modern mathematics by characterizing generalized (co)homology
theories like K-theories as a kind of knowledge-producing ``machines".
Probably most mathematicians would agree in that these machines had so
many useful applications that it seems a bit unfair to describe them
as mere chimeras.}

But for Connes such an ``ecumenical'' option is not available.  This
leads him into difficulties.  On the one hand, he relies upon the
Solovay model where all sets of real numbers are Lebesgue measurable
(see Subsection~\ref{what}), so as to relegate non-standard numbers to
the chimerical realm of mere \emph{tools}:
\begin{quote}
tout r\'eel non standard d\'etermine canoniquement un sous-ensemble
non Lebesgue mesurable de l'intervalle~$[0,1]$ de sorte qu'il est
impossible [Ste] d'en exhiber un seul (Connes 1997
\cite[p.~211]{Co97}).
\end{quote}
Here the reference ``[Ste]'' cited by Connes is an article by
J.~Stern~\cite{St}.  The main subject of Stern's article is a result
of S.~Shelah~(1984 \cite{Sh84}).  Shelah proved that the assumption of
the consistency of the proposition that all sets of real numbers are
Lebesgue measurable implies the consistency of inaccessible cardinals.
Connes' citation of Stern indicates that Connes was aware of Shelah's
1984 result.

On the other hand, Connes ignores the fact that for the consistency of
the proposition that all sets of real numbers are Lebesgue measurable,
Solovay (see Theorem~\ref{solt}) had to assume the existence of
inaccessible cardinals, and S.~Shelah showed that one cannot remove
the hypothesis of inaccessible cardinal from Solovay's theorem.
Meanwhile, Connes' meta-mathematical speculations, such as the claim
that ``noone will ever be able to name, etc.'' (see
Subsection~\ref{book}) rely on Solovay's theorem.  Therefore
ultimately Connes' meta-mathematical speculations rely on inaccessible
cardinals, as well.  The linchpin that keeps Connesian Platonism from
unraveling turns out to be an inaccessible cardinal, yet another
chimera.

What kind of evidence does Connes present in favor of his approach?
It is of two kinds:
\begin{enumerate}
\item
G\"odel's incompleteness theorem and Goodstein's theorem;
\item
feelings of eternity.
\end{enumerate}
We will examine these respectively in Subsections~\ref{103} and
\ref{104}.

\subsection{The incompleteness theorem: evidence for Platonism?}
\label{103}

\begin{figure}
\includegraphics[height=1.7in]{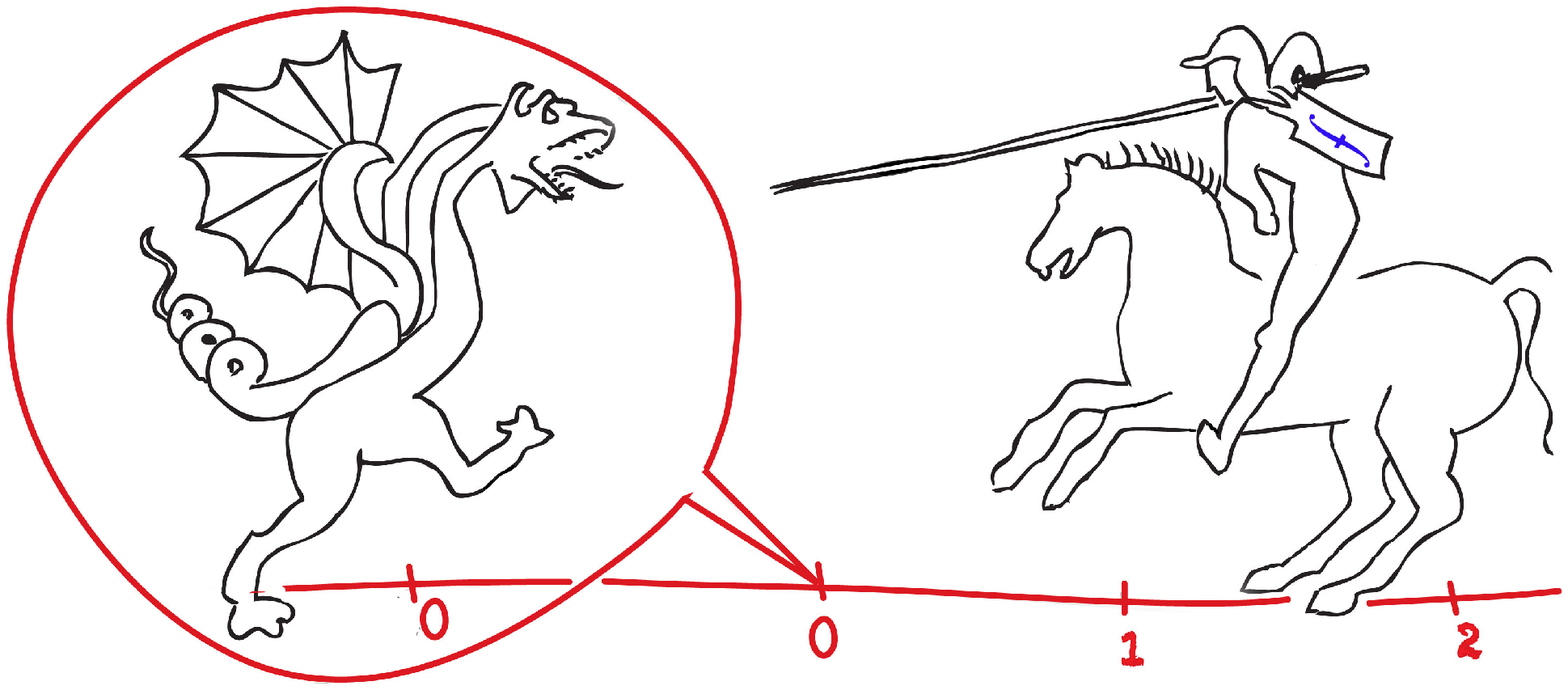}
\caption{\textsf{A virtual view of primordial mathematical reality: An
attempted slaying of a hyperreal chimera, following P.~Uccello}}
\label{uccello}
\end{figure}

There is an instance of apparent circular reasoning in one of Connes'
arguments in favor of his philosophical approach in the \emph{La
Recherche} interview \cite{Co00d}.%
\footnote{The discussion in this subsection was inspired by
I.~Hacking's \emph{The Mathematical Animal} \cite[chapter~5]{Ha13}.}
More specifically, Connes claims that G\"odel's incompleteness theorem
furnishes evidence in favor of Connes' philosophical approach, in that
it asserts the existence of ``true" propositions about natural numbers
that cannot be proved:
\begin{quote}
Or le th\'eor\`eme de G\"odel est bien plus m\'echant que cela.  Il
dit qu'il y a aura toujours une proposition vraie qui ne sera pas
d\'emontrable dans le syst\`eme.  Ce qui est beaucoup plus
d\'erangeant (Connes 2000 \cite{Co00d}).
\end{quote}
Such ``true'' propositions, undecidable in Peano Arithmetic (PA), are
taken by Connes to furnish evidence in favor of the hypothesis of a
mind-independent (Platonic) primordial mathematical reality (PMR),
referred to as \emph{r\'ealit\'e math\'ematique archa\"\i que} in the
interview.%
\footnote{An attempt to illustrate this concept graphically may be
found in Figure~\ref{uccello}, and further discussion in
Subsection~\ref{chimera}.}

However, the ``truth" of such propositions refers to truth relative to
an \emph{intended interpretation} of natural numbers, such as the one
built in Zermelo-Fraenkel set theory (ZF) or a fragment~$\text{ZF}_0$
thereof.  Relative to such an interpretation, the said propositions
are ``true" but not provable in PA.  At variance with Connes, K.~Kunen
presents G\"odel's theorem (in the context of ZF) in a philosophically
neutral way as follows:

\begin{quote}
if~$T$ is any consistent set of axioms extending ZF,%
\footnote{Actually it is sufficient to assume that~$T$ is consistent
and contains a suitable small set of axioms governing addition and
multiplication of natural numbers.}
then [the set]~$\{\varphi:T\vdash\varphi\}$ is not recursive \ldots{}
A consequence of this is G\"odel's First Incompleteness
Theorem---namely, that if such a~$T$ is recursive, then it is
incomplete in the sense that there is a sentence~$\varphi$ such that
$T\not\hskip.8pt\vdash\varphi$ and~$T\not\hskip.8pt\vdash\neg\varphi$
(Kunen 1980 \cite[p.~38]{Ku}).%
\footnote{The string~$T\vdash\varphi$ denotes the statement
``sentence~$\varphi$ is provable in theory~$T$'', while the
string~$T\not\hskip.8pt\vdash\varphi$ denotes the
statement~``$\varphi$ is \emph{not} provable in~$T$''.}

\end{quote}
With regard to Platonism, Kunen specifically mentions that G\"odel's
theorem, as well as the closely related Tarski's theorem on
non-definabi\-lity of truth, admit of platonist \emph{interpretations}
(rather than furnishing \emph{evidence} in favor of Platonism):
\begin{quote}
The platonistic interpretation of [Tarski's theorem] is that no
formula~$\chi(x)$ can say ``$x$ is a true sentence''%
\footnote{In some rare cases, it is possible to document a kind of
``model-theoretic failure" of the Tarski truth undefinability theorem.
Thus Kanovei \cite{Ka80} (and independently L.~Harrington,
unpublished) showed that in a suitable model of ZFC, the set of all
analytically definable reals is defined analytically; namely, it is
equal to the set of G\"odel-constructible reals.}
(Kunen 1980 \cite[p.~41]{Ku}).
\end{quote}

While Connes' argument appears to rely on an unspoken hypothesis of an
imbedding of such a fragment~$\text{ZF}_0$ in his PMR, he is certainly
free to believe in the hypothesis of such an imbedding
\begin{equation}
\label{101}
\text{ZF}_0\hookrightarrow\text{PMR}.
\end{equation}
Our goal here is to argue neither in favor nor against Connes'
hypothesis~\eqref{101}, but rather to point out an apparent
circularity inherent in Connes' argument.  Connes seeks to argue in
favor of Platonism based on G\"odel's result, but an unspoken
hypothesis of his argument is\ldots{} Platonism itself, about some
fragment~$\text{ZF}_0$ properly containing PA, betraying an apparent
circularity in his logic.

When Postel-Vinay (the \emph{La Recherche} interviewer) pressed Connes
for examples of statements that are ``true" but not provable, Connes
fell back on what he called ``La fable du li\`evre et de la tortue''
(``the hare and the turtle" phenomenon).  What Connes describes here
is in fact Goodstein's theorem \cite{Goo}.  As its name suggests, this
``true" theorem does admit of a proof, namely Goodstein's.  The proof
takes place not in PA but rather in a fragment
assuming~$\epsilon_0$-transfinite induction.  Relative to such a
widely accepted infinitary hypothesis, Goodstein's theorem is provable
and therefore true.

M.~Davis \cite{Da06} argued that~$\Pi^0_1$ sentences such as Cons(PA)
are equivalent to checking specific Diophantine questions, and
therefore their truth value should be determinate, and described such
a viewpoint as \emph{pragmatic Platonism} \cite{Da12b}.  Meanwhile,
Connes is characteristically evasive as to the scope of his platonist
beliefs, but his categorical tone suggests a rather broad Platonism.
What is clear, at any rate, is that his Platonism transcends
the~$\Pi^0_1$ class of the arithmetic hierarchy (since Goodstein's
theorem falls outside that class) and is probably much broader.  In
terms of Shapiro's distinction between \emph{realism in ontology} and
\emph{realism in truth-value} \cite[p.~37]{Sha}, Davis may be
described as a truth-value realist while Connes, an ontological one.

\subsection{Premonitions of eternity}
\label{104}

Connes' additional argument invokes ``a feeling of eternity" in
connection with his PMR:

\begin{quote}
La diff\'erence essentielle \ldots{} c'est qu'elle \'echappe \`a toute
forme de localisation dans l'espace ou dans le temps.  Si bien que
lorsqu'on en d\'evoile ne serait-ce qu'une infime partie, on \'eprouve
un \emph{sentiment d'\'eternit\'e}.  Tous les math\'ematiciens le
savent (Connes 2000 \cite{Co00d}) [emphasis added--the authors].
\end{quote}

Taking such a \emph{``sentiment d'eternit\'e''} as the ultimate litmus
test for one's reflection on what mathematics is and what
mathematicians do is a powerful means of effectively cutting off any
further reflection on the nature and the aim of mathematics and its
role in the context of culture and society at large.  After all, my
``sentiments" may be different from yours, and there is no room for
rational argumentation.  To take this road, one must invoke other
means of deciding which sentiments are justified and which are not,
such as appeals to the \emph{great mathematicians}: \emph{their}
``sentiments" are taken to need no justification at all, as they are
the only ones taken to have a legitimate say on what mathematics in
its essence really is (see, however, Subsection~\ref{atiyah} for the
anti-Platonist sentiments of M.~Atiyah).

However, relying on ``sentiments" when dealing with ontological issues
concerning mathematics not only has damaging effects on the discourse
about mathematics in general.  It also affects rather concrete issues
concerning the history of mathematics.  Arguably, a brand of
prescriptive Platonism about the real number continuum may, in fact,
be at the root of historical misconceptions concerning key figures and
pivotal mathematical developments.  Thus, consider the issue of
Fermat's technique of \emph{adequality} (stemming from Diophantus's
\parisotes) for solving problems of tangents and maxima and minima.
Fermat's technique involves an aspect of approximation and
``smallness'' in an essential way, as shown by its applications to
transcendental curves and variational problems such as Snell's law
(see Cifoletti~\cite{Ci}; Katz, Schaps, \& Shnider~\cite{KSS13}).
This aspect of Fermat's technique is, however, oddly denied by such
Fermat scholars as H.~Breger \cite{Bre94} and K.~Barner~\cite{Bar}.
Similarly, the non-Archimedean nature of Leibniz's infinitesimals is
routinely denied by some modern scholars (see Ishiguro \cite{Is90},
Levey \cite{Le08}), inspite of ample evidence is Leibniz's writings
(see Jesseph \cite{Je11}; Katz \& Sherry \cite{KS1}, \cite{KS2}).  A
close textual analysis of Cauchy's foundational writings reveals the
existence of a Cauchy--Weierstrass \emph{dis}continuity rather than
continuity, \emph{pace} Grabiner \cite{Gr81} (see B\l aszczyk et
al.~\cite{BKS}; Borovik \& Katz \cite{BK}; Br\aa ting \cite{Br}; Katz
\& Katz \cite{KK11a}, \cite{KK11b}; Sinaceur \cite{Si}).

\subsection{Cantor's dichotomy}

Cantor may be said to have opened Pandora's box of the ``chimeras" of
modern mathematics.  It appears that Cantor had a more elaborate and
flexible concept of mathematical reality than does Connes.  In his
\emph{Foundations of a general theory of manifolds} \cite{Can}, Cantor
pointed out that we may speak in two distinct ways of the reality or
existence of mathematical concepts.

First, we may consider mathematical concepts as \emph{real} insofar as
they, due to their definitions, occupy a fully determined place in our
mind whereby they can be distinguished perfectly from all other
components of our thought to which they stand in certain relations.
Thereby they are real since they may modify the substance of our mind
in certain ways.  Cantor called this kind of mathematical reality
\emph{intrasubjective} or \emph{immanent reality}.

On the other hand, one may ascribe reality to mathematical concepts
insofar as they can be considered as expressions or images of
processes and relations of the outside world.  Cantor referred to this
kind of reality as \emph{transient reality}.  Cantor had no doubt that
these two kinds of reality eventually came together.  Namely, concepts
with solely \emph{immanent reality} would, in the course of time,
acquire \emph{transient reality}, as well.  By this two-tiered concept
of the reality or ``Wirklichkeit" of mathematical entities Cantor
thought to have done justice to the idealist as well as to the realist
aspects of mathematics and mathematized sciences.%
\footnote{Cantor's actions did not always faithfully reflect his
professed flexible and tolerant attitude toward immanent ``chimeras".
As is well-known, he was eagerly hunting down infinitesimals of all
kinds as allegedly noxious chimeras to be eliminated.  One of his
strategies of elimination was the publication of a ``proof'' of an
alleged inconsistency of infinitesimals.  Accepting Cantor's analysis
on faith, B.~Russell declared infinitesimals to be inconsistent
\cite[p.~345]{Ru03}, influencing countless other philosophers and
mathematicians.  The errors in Cantor's ``proof'' are analyzed by
P.~Ehrlich~\cite{Eh06}.  It is interesting to note that Cantor's
contemporary B.~Kerry was apparently unconvinced by either Cantor's
feelings of eternity or by his ``proof'', and tried to put up an
argument, but was scornfully rebuffed by Cantor, who condemned Kerry's
alleged ``deplorable psychologistic blindness'' (see C.~Proietti
\cite[p.~356]{Pr}) and concluded: ``\emph{Dixi et salvavi animam
meam}.  I think I did my best to dissuade you from your deplorable
mistakes'' (ibid.).}

Our analysis of Connes' approach should not be misunderstood.  We do
not deny that the distinction between \emph{tool} and \emph{object} is
an eminently useful one.  The point is that one has to take into
account the historical and relative character of this distinction.
Exactly this Connes' Platonism does not do.  Thereby it is blinded to
certain essential features of modern mathematical knowledge.  The
manifest historical evolution of the domain of mathematical objects
and the emergence of new tools, which depend on the \emph{changing}
character of the object domain, points to a dynamism of the
ontological realm of mathematics to which Connes' vision of a
``primordial mathematical reality" (PMR) is directly opposed.  Connes'
account of mathematical knowledge implies a static ontology.  The
innate weakness of Connes' vision of PMR is that it ignores the
inevitable interaction between tools and objects in science.

Furthermore, such an interaction between tools and objects brings into
play the institution of a \emph{subject} that is actively using and
creating both tools and objects for its specific purposes that may
change over time and historical context.  In Connes' account, the
subject (that is engaged in ``doing" mathematics) fatally resembles
the ideal, non-empirical subject of classical philosophy for which
finiteness and other empirical limitations of the real empirical
subjects were philosophically irrelevant.

Despite his platonist preferences, history as well as
subject-with-a-history is surreptiously introduced by Connes himself,
however.  The talk of \emph{tools} only makes sense if a
\emph{subject}, i.e., an agent is presupposed that \emph{employs}
these tools for its purposes.  Connes' subject is a transmundane and
very abstract entity.
%
%
A more convincing choice of the subject would be a historically
situated subject.  After all, it can hardly be denied that mathematics
as every other scientific discipline has undergone a historical
development; \emph{our} mathematics is not the same as Greek
mathematics, and it is hardly plausible that the mathematics of the
future will be ``essentially the same" as present-day mathematics.
The line between \emph{tools} and \emph{objects} is moving.  A tool
may gain the status of an object and, conversely, an object may become
a tool in a suitable context.

\subsection{Atiyah's anti-platonist realism}
\label{atiyah}

Not all great contemporary mathematicians share Connes' philosophical
position.  Thus, Sir Michael Atiyah confided:

\begin{quote}
I consider myself as a realist.  I think the mathematics we use is
derived from the outside world by observation and abstraction. If we
didn't live in the outside world and see things, we wouldn't have
invented things and thought of things as we do. I think that much of
what we do is based on what we see, but then abstracted and
simplified, and in that sense they become the ideal things of Plato,
but they have an origin in the outside world and that's what brings
them close to physics.  \ldots{} You can't separate the human mind
from the physical world. And therefore everything we think of, in some
sense or other, derives from the physical world (Atiyah 2005,
\cite[p.~38]{At}).
\end{quote}

Atiyah's outline of a realistic conception of mathematics is not, of
course, without problems.  For instance, one may object that we do not
spend our life time by merely ``seeing the outside world".  Rather, we
are beings in a material world and have to come to terms with the
multifarious challenges that the world poses to us.  Hence, rather
than describing our contact with the outside world as ``seeing", it
may be more appropriate to adopt a broader approach that emphasizes
the multifaceted totality of the various activities in which cognizing
beings like us are engaged.  One may object that Atiyah does not
elaborate much on the profound issue of what exactly is meant by
``deriving mathematics from the outside world" and how this is carried
out.  We think that such a criticism would be a bit unfair.  One may
well argue that these issues are not, properly speaking, mathematical
issues and therefore are not a primary concern for mathematicians.
%

\subsection{Mac Lane's form and function}

A more elaborate account of how ``mathematics is derived from the
outside world" can be found in Saunders Mac Lane's \emph{Mathematics,
Form and Function} (Mac Lane 1986 \cite{Mac}).  This book recorded Mac
Lane's
\begin{quote}
efforts \ldots{} to capture in words a description of the form and
function of Mathematics, as a background for the Philosophy of
Mathematics" (Preface).
\end{quote}
Here Mac Lane compiled a list of rather mundane activities such as
collecting, counting, comparing, observing, moving and others that can
be considered as the modest origins of the high-brow concepts of
contemporary mathematics (ibid., p.~35).  An interesting elaboration
of Mac Lane's account
%
%
may be found in \emph{Where Mathematics Comes From.  How the Embodied
Mind Brings Mathematics into Being} (Lakoff \& N\'u\~nez 2000
\cite{LN}).

The details of the processes underlying the historical evolution of
mathematics may not be fully understood yet.
%
%
However, cutting off any further discussion on these issues by falling
back on ``feelings of eternity" does not seem the best way to meet
such challenges.  Mathematics, as any other intellectual endeavor,
cannot be considered as an autonomous domain totally cut off from
other areas of knowledge.  As Atiyah put it explicitly:
\begin{quote}
The idea that there is a pure world of mathematical objects (and
perhaps other ideal objects) totally divorced from our experience,
which somehow exists by itself is obviously inherent nonsense (Atiyah
2006 \cite[p.~38]{At}).
\end{quote}
A PMR-free perspective on mathematics is gaining momentum.  In fact,
Connes' feelings of eternity may be misdirected.  Scholars from many a
discipline converge to a view that thinking about mathematics should
not treat the latter as an isolated endeavor, separate from other
areas of knowledge.

\subsection{Margenau and Dennett: To be or \ldots}

Connes' radical Platonism with its postulation of a strict separation
of the sphere of mathematics from the rest of the world is, in a
sense, radically anti-modern.  Modernity in the sciences began with a
turn toward epistemological and semantical questions, leaving aside
classical ontological questions such as ``What is the essence of the
world?", ``What is the essence of Man?", or, more to the point of the
present paper, ``What is the essence of number or space?".  Instead,
in the modern perspective, semantical and epistemological questions
such as ``What is the meaning of this or that scientific concept in
this or that context?", ``What is scientific knowledge?", or ``Can one
make sense of the progress in science?" take centerstage.  In this
way, ontology, epistemology, and semantics get inextricably
intertwined.  In particular, ontology became theory-dependent.  For
the mathematized sciences of nature, the neo-Kantian philosopher Ernst
Cassirer expressed this observation explicitly as follows:

\begin{quote}
[Scientific] concepts are valid not in that they copy a fixed, given
being, but insofar as they contain a plan for possible constructions
of unity, which must be progressively verified in practice \ldots{}
(Cassirer 1957 \cite[p.~476]{Ca57}).
\end{quote}

What we need is not the objectivity of absolute concepts (it seems
difficult to give convincing arguments to account for how one could
have cognitive access to such concepts), but rather objective methods
which determine the rational and reliable practice of our
intersubjective empirical science.  As Cassirer put it,

\begin{quote}
What we need is not the objectivity of absolute objects, but rather
the objective determinacy of the \emph{method of experience} (ibid.)%
%
%
\end{quote}

Cassirer's characterisation of scientific concepts as applied to
mathematical concepts amounts to the contention that mathematical
concepts should not be conceived of as intending to copy a
pre-existing platonic universe but ``contain plans for possible
constructions of unity".  This characterization would match quite well
with Hilbert's dictum ``By their fruits ye will known them".  If this
is true, a ``theory of chimeras" \`a la Connes hardly provides a
promising framework for dealing with these problems.

Rather, what is needed is an investigation of the entire spectrum of
the various meanings of the concept of \emph{being} as it is used in
modern science.  The need for such an investigation was pointed out by
Cassirer's friend and colleague, the renowned physicist Henry
Margenau, by means of the following provocative question:

\begin{quote}
Do masses, electrons, atoms, magnetic field strengths etc., exist?
Nothing is more surprising indeed than the fact that \ldots{} most of
us still expect an answer to this question in terms of yes or
no. \ldots{} Almost every term that has come under scientific scrutiny
has lost its initally absolute significance and acquired a range of
meaning of which even the boundaries are often variable.  Apparently
the word \emph{to be} has escaped this process (Margenau 1935
\cite[p.~164]{Ma35}).
\end{quote}

Margenau argued in favor of a nuanced concept of ``the real" based on
an elaborate theory of theoretical constructs in which ``tools" and
``objects" interact in complex ways (cf. Margenau 1935 \cite{Ma35},
Margenau 1950 \cite{Ma50}).

Sixty years later, Margenau's question was taken up and generalized to
the object of other sciences by
%
Daniel Dennett:

\begin{quote}
Are there really beliefs? Or are we learning (from neuro-science and
psychology, presumably) that strictly speaking, beliefs are figments
of our imagination, items in a superseded ontology. Philosophers
generally regard such ontological questions as admitting just two
possible answers: either beliefs exist or they do not. (Dennett 1991
\cite[p.~27]{De}).
\end{quote}

Dennett argued that an ontological account centered around the concept
of ``patterns" may be helpful to develop an ``intermediate''
(Dennett's term) position that conceives of beliefs and other
questionable abstract entities as patterns of some data.  Taking data
as a bit stream, a pattern is said to exist in some data, i.e., is
real if there is a description of the data that is more efficient than
the bit map, whether or not anyone can concoct it. Thereby centers of
gravity exist in physicalist ontologies because they are good abstract
concepts that perform some useful work.  Meanwhile, bogus concepts
such as ``Dennett's lost socks center" (defined as ``the center of the
smallest sphere that can be circumscribed around all the socks Dennett
ever lost in his life") do not obtain this status but remain
meaningless ``chimeras" (ibid., 28).

In a somewhat analogous way, Michael Resnik and other philosophers of
mathematics are working on a project of describing ``mathematics as a
science of patterns", in which Resnik defends the thesis that
mathematical structures obtain their reality as ``patterns of reality"
(Resnik 1997 \cite{Re}).

This section is not the place to engage in an in-depth study of these
and similar attempts to clarify the murky issue of the ontology and
epistemology of mathematics.  Our goal is merely to evoke some
possibly fruitful directions of inquiry that may help overcome the
limitations of the traditional accounts of formalism, intuitionism,
and platonism.  In the long run it is unsatisfying (to put it mildly)
to play off against each other these classical positions over and over
again, by manufacturing unappealing and unrealistic strawmen of the
other party.  Such dated ideas on the nature of mathematics do not
exhaust the spectrum of possible approaches to the epistemology and
ontology of mathematics.

Connes' views on non-standard analysis are inseparable from his
philosophical position, as we discuss in Section~\ref{two}.

\section{``Absolutely major flaw'' and ``irremediable defect''}
\label{two}

Having clarified the philosophical underpinnings of Connes' views in
Section~\ref{otte}, we now turn to the details of his critique.
Connes published his magnum opus \emph{Noncommutative geometry} (an
expanded English version of an earlier French text) in 1994.  Shortly
afterwards, Connes published his first criticism of non-standard
analysis (NSA) in 1995, describing the non-standard framework as being
``inadequate".  In 1997, the adjective was \emph{``d\'ecevante''} (see
\cite{Co97}).  By 2000, Connes was describing non-standard numbers as
``chimeras".  Such criticisms have appeared in his books, research
articles, interviews, and a blog.

It is instructive to compare two papers Connes wrote around 2000.  The
paper \cite{Co00c} in \emph{Journal of Mathematical Physics} (JMP)
presents Connes' theory of infinitesimals without a trace of any
reference to either NSA or the Solovay model.  The other text from the
same period (see \cite{Co99}, \cite{Co00b}, \cite{Co00}) presents the
-- by then -- familiar meta-mathematical speculations around the
Solovay model (see Section~\ref{definable} for details), and proceeds
to criticize NSA.  The JMP text demonstrates that Connes is perfectly
capable of presenting his approach to infinitesimals (which he claims
to be entirely different from Robinson's) without criticizing NSA.

Connes was familiar with the ultrapower
construction~$\R^\N/\mathcal{F}$ of the hyperreals, having authored
the 1970 articles \cite{Co70a} and \cite{Co70b}.  At least on one
occasion, Connes described ultraproducts as ``very efficient'',%
\footnote{See main text at footnote~\ref{uccello2}.}
which adds another dimension to the puzzle.  To understand Connes'
position, one may have to examine the historical context of his
changing attitude toward non-standard analysis.  After Robinson's
death in 1974, many voices were heard that were critical of Robinson's
theory.  Active in this area were Paul Halmos and his student Errett
Bishop \cite{Bi77} (see Katz \& Katz \cite{KK11d}).  Some of the
criticisms were plain incoherent, such as John Earman's \cite{Ea} in
1975 (see Katz \& Sherry \cite[Section~11.2]{KS1}), suggesting that
for a time, it was sufficient to criticize Robinson to get published.
It may have become difficult starting in the mid-1970s to be a
supporter of Robinson, and it would have been natural for young
researchers to seek to distance themselves from him.  The objection to
hyperreal numbers on the part of many mathematicians may be due,
consciously or unconsciously, to their attitude that the traditional
model of the real numbers in the context of ZF is a true
representation of Reality itself%
\footnote{An alternative view is explored in \cite{KK11c}.}
and that hyperreal numbers are therefore a contrived model that does
not represent anything of interest, even if it provides a solution to
some paradox.  E.~Nelson, however, turned the tables on this attitude,
by introducing an enriched syntax into ZF, building the ``usual" real
line~$\R$ in ZF with the enriched syntax, and exhibiting
infinitesimals within the real line~$\R$ itself (see Nelson
\cite{Ne}).  Related systems were elaborated by K.~Hrb\'a\v{c}ek
(1978, \cite{Hrb}), T.~Kawai (1983, \cite{Kaw}) and Kanovei (1991,
\cite{Ka91}).

\subsection{The book}
\label{book}

The 2001 book \cite{CLS} was ostensibly authored by Connes,
A.~Lichnerowicz, and M.~\Schutzenberger.  Lichnerowicz and
\Schutzenberger{} died several years prior to the book's publication.
A reviewer notes:
\begin{quote}
The main contributions to the conversations come from Connes [\ldots]
and the fact that some of Connes' contributions look relatively
polished may indicate that they have been edited to some extent
[\ldots] Connes often explains a topic in a more or less systematic
way; \Schutzenberger{} makes interesting comments, often from a very
different angle while introducing many side-subjects, Lichnerowicz
interjects skeptical remarks (D.~Dieks 2002 \cite{Die}).
\end{quote}
The book's discussion of NSA in the form of an exchange with
\Schutzenberger{} appears on pages 15-21.  Here Connes expresses
himself as follows on the subject of non-standard analysis:

\begin{quote}
\quad A.C. - [\ldots] I became aware of an absolutely major flaw in
this theory, an irremediable defect.  It is this: in nonstandard
analysis, one is supposed to manipulate infinitesimals; yet, if such
an infinitesimal is given, starting from any given nonstandard number,
a subset of the interval \emph{automatically} arises which is not
measurable in the sense of Lebesgue.

\medskip
M.P.S. - Aha!

\medskip
A.C. - Yes, a nonstandard number yields in a simple \emph{canonical}
way, a subset of~$[0,1]$ which is not measurable in the sense of
Lebesgue [\ldots] What conclusion can one draw about nonstandard
analysis?  This means that, since noone will ever be able to name a
nonstandard number, the theory remains \emph{virtual}, and has
absolutely no significance except as a tool to understand ``primordial
mathematical reality''%
\footnote{See Subsection~\ref{103} for an analysis of the term
``primordial mathematical reality''.}
(Connes 2001, \cite[p.~16]{CLS}) [emphasis added--the authors]
\end{quote}
Connes goes on%
\footnote{The continuation of the discussion is dealt with in
Subsection~\ref{contrast}.}
to criticize the role of the axiom of choice in non-standard analysis
(ibid., p.~17).

Connes' criticisms of non-standard analysis have appeared in numerous
venues, and have been repeatedly discussed.%
\footnote{See, e.g.,
http://mathoverflow.net/questions/57072/a-remark-of-connes}
Some of the epithets he used for NSA, arranged by year, appear in
Table~\ref{epithets}.

\renewcommand{\arraystretch}{1.3}
\begin{table}
\[
\begin{tabular}[t]
{ | p{.33in} || p{3in} | p{.9in} | p{.5in} | p{.75in} | p{.5in} |}
\hline date & epithet & source \\ \hline\hline 1995 & ``inadequate'' &
\cite[p.~6207]{Co95} \\ \hline 1997 & ``\emph{d\'ecevante}''
[disappointing] & \cite[p.~211]{Co97} \\ \hline 2000 & ``very bad
obstruction'' & \cite[p.~20]{Co00} \\ \hline 2000 & ``chimera'' &
\cite[p.~21]{Co00} \\ \hline 2001 & ``absolutely major flaw'';
``irremediable~defect''; ``the theory remains virtual'' &
\cite[p.~16]{CLS} \\ \hline 2007 & ``I have found a catch in the
theory''; ``it seemed utterly doomed to failure to try to use
non-standard analysis to do physics'' & \cite[p.~26]{Co07a} \\ \hline
2007 & ``the promised land for `infinitesimals'\,''; ``the end of the
rope for being `explicit'\,'' & \cite{Co07b} \\ \hline
\end{tabular}
\]
\caption{\textsf{Connes' epithets for NSA arranged chronologically.}}
\label{epithets}
\end{table}
\renewcommand{\arraystretch}{1}

Some of Connes' criticisms are more specific than others.  Thus, the
precise meaning of his terms such as ``virtual theory'' and
``primordial mathematical reality'' is open to discussion (see
Section~\ref{otte}).  We will focus on the more mathematically
identifiable claim of a \emph{canonical} derivation of a Lebesgue
nonmeasurable set from a non-standard number, as well as the role of
Solovay's models in Connes' criticism.

Note that a construction of a nonmeasurable set starting from a
hyperinteger was described decades earlier by W.~Luxemburg (1963
\cite{Lu63} and 1973 \cite[Theorem~10.2, p.~66]{Lu73}), K.~Stroyan \&
Luxemburg (1973 \cite{SL}), and M.~Davis (1977
\cite[pp.~71-74]{Da77}).%
\footnote{Davis noted recently \cite{Da12} that he based his
construction on (Luxemburg 1962~\cite{Lu62}), by filling in the proof
of Theorem~9.1 in \cite[p.~72]{Da77} and otherwise following
Luxemburg.}

\subsection{Skolem's non-standard integers}

Before going into the \emph{mathematical} details of Connes' criticism
of non-standard numbers, we would like to comment on its
\emph{historical} scope.  Connes' criticism of non-standard integers
is worded in such a general fashion that one wonders if it would
encompass also the non-standard integers constructed by T.~Skolem in
the 1930s (see Skolem 1933 \cite{Sk33}, 1934 \cite{Sk34}; an English
version may be found in Skolem 1955 \cite{Sk55}).  Skolem's
accomplishment is generally regarded as a major milestone in the
development of 20th century logic.

D.~Scott \cite[p.~245]{Sc} compares Skolem's predicative approach with
the ultrapower approach (Skolem's nonstandard integers are also
discussed by Bell \& Slomson \cite{BS} and Stillwell
\cite[pp.~148-150]{St77}).  Scott notes that Skolem used the ring~$DF$
of algebraically (first-order) definable functions from integers to
integers.  The quotient~$DF/P$ of~$DF$ by a minimal prime ideal~$P$
produces Skolem's non-standard integers.  The ideal~$P$ corresponds to
a prime ideal in the Boolean algebra of idempotents.  Note that the
idempotents of~$DF$ are the characteristic functions of (first-order)
definable sets of integers.  Such sets give rise to a denumerable
Boolean algebra~$\mathfrak{P}$ and therefore can be given an
\emph{ordered basis}.  Such a basis for~$\mathfrak{P}$ is a nested
sequence%
\footnote{\label{f11}We reversed the inclusions as given in
\cite[p.~245]{Sc} so as to insist on the analogy with a filter.}
\[
X_n\supset X_{n+1} \supset \ldots
\]
such that~$Y\in\mathfrak{P}$ if and only if~$Y\supset X_n$ for a
suitable~$n$.  Choose a sequence~$(s_n)$ such that
\[
s_n\in X_{n}\setminus X_{n+1}.
\]
Then functions~$f,g\in DF$ are in the same equivalence class if and
only if 
\[
(\exists N) (\forall n\geq N)\, f(s_n)=g(s_n).
\]
The sequence~$(s_n)$ is the \emph{comparing function} used by Skolem
to partition the definable functions into congruence classes.  Note
that, even though Skolem places himself in a context limited to
definable functions, a key role in the theory is played by the
comparing function which is \emph{not} definable.

\subsection{The Connes character}

In Subsection~\ref{book}, we cited Connes to the effect that a
nonmeasurable set ``automatically'' arises, and that a non-standard
number ``canonically'' produces such a set.  Challenged to elaborate
on his claim, Connes expressed himself as follows:

\begin{quote}
Pour exhiber un ensemble non-mesurable a partir d'un entier
non-standard~$n$ il suffit de prendre le caract\`ere
de~$G=(\Z/2\Z)^\N$ qui est donn\'e par l'evaluation de la
composante~$a_n$ \ldots On obtient un caract\`ere non continu de~$G$
et il est donc non-mesurable (Connes 2009, \cite{Co09}).
\end{quote}
Similar remarks appear at Connes' non-standard blog (Connes
\cite{Co07b}).

In more detail, consider the natural numbers~$\N$, and form the
infinite product~$G=(\Z/2\Z)^\N$ (when equipped with the product
topology, it is homeomorphic to the Cantor set).
%
%
Each~$n\in\N$ gives rise to a homomorphism~$\chi_n:G\to\Z/2\Z$ given
by evaluation at the~$n$-th component.  Each element~$x\in G$ can be
thought of as a map
\begin{equation}
\label{31c}
x: \N \to \Z/2\Z=\{e,a\},
\end{equation}
where~$e$ is the additive identity element and~$a$ is the
multiplicative identity element.  Consider the set
\begin{equation}
\label{32c}
A_x=x^{-1}(a)\subset\N.
\end{equation}
Then~$x$ can be thought of as an ``indicator'' function of the
set~$A_x$.  In non-standard analysis, the map~$x$ of~\eqref{31c} has a
natural extension~${}^*x$ whose domain is the ring of hypernatural
numbers,~${}^*\N$:
\begin{equation}
\label{32b}
{}^*x: {}^*\N \to \Z/2\Z.
\end{equation}
Now let~$n\in{}^*\N\,\setminus\N$ be an infinite hypernatural.  The
evaluation of the map~${}^*x$ of~\eqref{32b} at~$n$ gives the
value~${}^*x(n)\in\Z/2\Z$ of~${}^*x$ at~$n$.  This again produces a
homomorphism from~${}^*G$ to~$\Z/2\Z$.  Its restriction to~$G\subset
{}^*G$ is denoted
\begin{equation}
\label{32}
\chi_n:G\to\Z/2\Z,\quad x\mapsto{}^*x(n).
\end{equation}
Thus, the character%
\footnote{A \emph{character} is generally understood to have image
in~$\C$; if one wishes to think of~\eqref{32} as a character, one
identifies~$\Z/2\Z$ with~$\{\pm 1\}\subset\C$.}
$\chi_n$ maps~$G$ to~$\Z/2\Z=\{e,a\}$.  Here
\begin{equation}
\label{31b}
\chi_n(x)=a \text{ if and only if } n\in {}^*\!A_x,
\end{equation}
where~${}^*\!A_x \subset{}^*\N$ is the natural extension of the
set~$A_x\subset\N$ of~\eqref{32c}.  Connes notes that the
character~$\chi_n$ is nonmeasurable.
%
%
He describes the passage from~$n$ to the character as ``canonical'',
and alleges that non-standard analysis introduces entities that lead
``canonically'' to nonmeasurable objects.%
\footnote{\label{f9}Another interpretation:~$G=(\Z/2\Z)^\N$ is the
standard product which is a compact metrizable group.  Each
element~$x\in G$ has an internal extension~${}^*\!x$ defined
on~${}^*\N$.  Thus, if~$n$ is a standard or non-standard hypernatural,
then~${}^*\!x$ can be evaluated at~$n$.  Now the continuous dual
of~$G$, by Pontryagin duality, is the algebraic direct sum of
countably many copies of~$\Z/2\Z$ with the discrete topology.  Thus,
the evaluation at a non-standard integer~$n$ is not continuous and
therefore not measurable,
%
%
and cannot be equal a.e.~to a Borel function.}

\subsection{From character to ultrafilter}
\label{four}

The Connes character~$\chi_n^{\phantom{I}}$ carries the same
information as an ultrafilter.  Indeed, consider the inverse image
of~$a\in\Z/2\Z$ under the character~$\chi_n^{\phantom{I}}$
of~\eqref{32}, namely,~$\chi_n^{-1}(a)\subset G$.  To
each~$x\in\chi_n^{-1}(a)$, we can associate the
subset~$A_x\in\mathcal{P}(\N)$ of~\eqref{32c}.%
\footnote{Here~$\mathcal{P}(\N)$ denotes the set of subsets of~$\N$.}
If~$n\in{}^*\N\setminus\N$ is a fixed hypernatural, then the
collection
\[
\left\{A_x \in \mathcal{P}(\N) : \chi_{n}^{\phantom{I}}(x)=a \right\}
\]
yields a free ultrafilter on~$\N$.  By~\eqref{31b}, Connes'
construction can be canonically identified with the following
construction.

\begin{construction}
\label{21}
Choose an unlimited hypernatural~$n\in{}^*\N$, and construct the
ultrafilter~$\mathcal{F}\subset\mathcal{P}(\N)$ consisting of
subsets~$A\subset\N$ whose natural extension~${}^*\!A\subset{}^*\N$
contains~$n$:
\begin{equation}
\label{41}
\mathcal{F}= \left\{ A\in \mathcal{P}(\N) : n\in {}^*\!A \right\}.
\end{equation}
\end{construction}

The important remark at this stage is that Connes' construction
exploits a new principle of reasoning introduced by Robinson, called
the \emph{transfer principle}.%
\footnote{The transfer principle for ultraproduct-type nonstandard
models follows from \L o\'s's theorem dating from 1955 (see
\cite{Lo55}).}
The reliance of the construction on the transfer principle was
acknowledged%
\footnote{\label{tp}See Subsection~\ref{contrast} at footnote~\ref{f7}
for a further discussion of the role of the transfer principle.}
by Connes \cite{Co12}.

\begin{remark}
\label{42b}
If one applies Construction~\ref{21} to the hypernatural
\begin{equation}
\label{31}
n=[(1,2,3,...)],
\end{equation}
i.e., the equivalence class of the sequence listing all the natural
numbers, then one recovers precisely the ultrafilter~$\mathcal{F}$
used in the ultrapower construction of a hyperreal field as the
quotient%
\footnote{More precisely, we form the quotient of~$\R^\N$ by the space
of real sequences that vanish on members of~$\mathcal{F}$.  The
notational ambiguity is widespread in the literature.}
\begin{equation}
\label{42}
\IR=\R^\N/\mathcal{F}.
\end{equation}
\end{remark}

\subsection{A forgetful functor}
\label{five}

Connes has repeatedly used the terminology of ``canonical'' in his
publications, as in the claim that ``a hyperreal number
\emph{canonically} produces'' a nonmeasurable entity.  To an
uninformed reader, this may sound similar to an assertion that ``to
every rational number one can \emph{canonically} associate a pair of
integers" (reduce to lowest terms), or ``to every real number one can
\emph{canonically} associate a unique Dedekind cut" on~$\Q$.  Both of
these statements are true if the field is given up to isomorphism,
with no additional structure.

It is not entirely clear if Connes means to choose an element from a
specific model of a hyperreal field, or an element%
\footnote{\label{orbit}More precisely, the orbit of an element under
field automorphisms.}
from an isomorphism type of such a model (i.e., its class up to
isomorphism).  We will therefore examine both possibilities:
\begin{enumerate}
\item
element of an isomorphism type of a hyperreal field; or
\item
element of a particular non-standard model.
\end{enumerate}
Briefly, we argue that in the former case Connes' claim is false.
Meanwhile, in the latter case, the complaint is moot as we already
have an ultrafilter~$\mathcal{F}$, namely the one used to build the
model as in \eqref{42}.  Thus, Connes' ``canonical'' procedure is
canonically equivalent to a black box%
\footnote{\label{black1}See also main text in Section~\ref{virtual} at
footnote~\ref{black2}.}
that canonically returns its input (namely, the original
ultrafilter~$\mathcal{F}$; see Remark~\ref{42b}).  More precisely, it
is a forgetful functor~$\Phi$ from the category~$\mathcal{E}$ of
hyperreal enlargements to the category~$\mathcal{U}$ of ultrafilters:
\begin{equation}
\label{51}
\Phi: \mathcal{E}\to\mathcal{U}, \quad \Phi \left( \R; \mathcal{F};
\IR=\R^\N/\mathcal{F}; * \right) = \mathcal{F}.
\end{equation}

\subsection{$P$-points and Continuum Hypothesis}
\label{52}

We argue that to produce a canonical ultrafilter from a hyperreal, an
isomorphism type of~${}^*\R$ does not suffice.  To see this, assume
for the sake of simplicity the truth of the continuum hypothesis (CH);
note that a procedure claimed to be ``canonical" should certainly work
in the assumption of CH, as well.  Now in the traditional
Zermelo--Fraenkel set theory with the Axiom of Choice (ZFC) together
with the assumption of CH, we have the following theorem (see Erd\"os
et al.~1955, \cite{EGH}).

\begin{theorem}[Erd\"os et al.]
\label{erdos}
In ZFC+CH, all models of~$\IR$ of the form~$\R^\N/\mathcal{F}$ are
isomorphic as ordered fields.%
\footnote{In fact, the uniqueness up to isomorphism of this ordered
field is equivalent to CH (see Farah \& Shelah \cite{FS}).}
\end{theorem}

Meanwhile, the ultrafilter~$\mathcal{F}$ may or may not be of a type
called a ``$P$-point''.  The most relevant property of an
ultrafilter~$\mathcal{F}$ of this type is that every infinitesimal
in~$\R^\N/\mathcal{F}$ is representable by a null sequence, i.e., a
sequence tending to zero (see Cutland et al.~\cite{CKKR}).  Meanwhile,
not all ultrafilters are~$P$-points.%
\footnote{\label{rudin}Thus, W.~Rudin (1956, \cite{Ru}) proved the
following results assuming CH.  Recall that a space is called
homogeneous if for any two points, there is a homeomorphism taking one
to the other.  

Theorem~4.4:~$\beta\N-\N$ is not homogeneous;
Theorem~4.2:~$\beta\N-\N$ has~$2^c$~$P$-points; Theorem 4.7: for any
two~$P$-points of~$\beta\N-\N$, there is a homeomorphism
of~$\beta\N-\N$ that carries one to the other.}

Thus, the isomorphism type of~$\IR$ does not retain the information as
to which ultrafilter was used in the construction thereof.
If~$\mathcal{F}$ is a~$P$-point, then the hypernatural~\eqref{31} fed
into~\eqref{41} will return the~$P$-point ultrafilter~$\mathcal{F}$
itself, but also \emph{every} choice of a
hyperinteger~$n\in{}^{*}\N\setminus\N$ would yield a~$P$-point (this
follows from the properties of the Rudin--Keisler order on the
ultrafilters).%
%
%

If a~$P$-point~$\mathcal{F}$ were used in the construction
of~$\IR$, any imaginable ``canonical'' construction (such as
Connes', exploiting the transfer principle) would have to yield
a~$P$-point, as well.  But if all one knows is the isomorphism class
of~$\IR$, the nature of the ultrafilter used in the construction
cannot be detected; it may well have been a \emph{non}-$P$-point
ultrafilter.  We thus obtain the following:

\begin{quote}
\emph{There does not exist a canonical construction of a nonprincipal
ultrafilter from an element%
\footnote{See footnote~\ref{orbit}.}
in an isomorphism type of a hyperreal field.}
\end{quote}

Such a construction could not exist unless one is working with
\emph{additional} data (i.e., in addition to the isomorphism type),
such as a specific enlargement~$\R\to\IR$ with a transfer
principle, where we can apply Construction~\ref{21}.  However, the
construction of such an enlargement requires an ultrafilter to begin
with!  This reveals a circularity in Connes' claim.

\subsection{Contrasting infinitesimals}
\label{contrast}

We continue our analysis of the discussion between \Schutzenberger{}
and Connes started above in Subsection~\ref{book}.  Connes contrasts
\emph{his} infinitesimals with Robinson's infinitesimals in the
following terms:
\begin{quote}
An infinitesimal [in Connes' theory] is a certain type of operator
which I am not going to define.  What I want to emphasize is that in
the critique of the nonstandard model, the axiom of choice plays an
extremely important role that I would like to make explicit.  In
logic, when one constructs a nonstandard model, for example of the
integers, or of the real line, one tacitly uses the axiom of choice.
It is applied in an \emph{uncountable} situation (Connes
\cite[p.~17]{CLS}) [emphasis added--the authors].
\end{quote}
The comment appears to suggest that Robinson's theory relies on
uncountable choice but Connes' does not.  The validity or otherwise of
this suggestion will be discussed below (see end of this Subsection).
The discussion continues as folows:
\begin{quote}
\quad M.P.S. - What you are saying is fantastic.  I had never paid
attention to the fact that the countable axiom of choice differed from
the uncountable one.  I must say that I have nothing to do with the
axiom of choice in daily life.

\medskip

A.L. - Of course not! (Connes \cite[p.~17]{CLS}).
\end{quote}

What emerges from \Schutzenberger's comments is that he ``never paid
attention'' to the distinction between the countable case of the axiom
of choice and the general case.  The continuation of the discussion
reveals that \Schutzenberger{} is similarly ignorant of the concept of
a \emph{well ordering}:
\begin{quote}

\quad M.P.S. What do you mean by ``well ordering''!?

\medskip

A.C.  Well ordering!  The integers have the property that \ldots
[there follows a page-long introduction to well ordering.]

\medskip

M.P.S. Amazing!

\medskip

A.L. [Lichnerowicz] So the countable and uncountable axioms of choice
are different.

\medskip

A.C. Absolutely.  It is worth noting that most mainstream mathematics
only requires the countable axiom of choice%
\footnote{It is difficult to argue with a contention that ``mainstream
mathematics only requires the countable axiom of choice'', since the
term \emph{mainstream mathematics} is sufficiently vague to accomodate
a suitable interpretation with respect to which the contention will
become accurate.  Note, however, that such an interpretation would
have to relegate Connes' work in functional analysis on the
classification of factors (for which Connes received his Fields medal)
to the complement of ``most mainstream mathematics'', as his work
exploited ultrafilters in an essential manner, whereas ZF+DC is not
powerful enough to prove the existence of ultrafilters (see
Remark~\ref{42c}).}
[\ldots] (Connes \cite[p.~20-21]{CLS}).
\end{quote}

Connes' discussion of the distinction between the countable axiom of
choice (AC) and the general AC appears to suggest that one of the
shortcomings of non-standard analysis is the reliance on the
uncountable axiom of choice.

Such a suggestion is surprising, since Connes' own framework similarly
exploits nonprincipal ultrafilters which cannot be obtained with
merely the countable AC (see Remark~\ref{dc}, Section~\ref{dix}, and
Remark~\ref{von}).  The impression created by the discussion that
Connes' theory relies on countable AC alone, is therefore spurious.

\subsection{A virtual discussion}
\label{virtual}

Sh\"utzenberger was not in a position to challenge any of Connes'
claims due to ignorance of basic concepts of set theory such as the
notion of a well ordering.  Had he been more knowledgeable about such
subjects, the discussion may have gone rather differently.

\begin{quote}

\quad M.P.S. - I have the following question concerning the evaluation
at a nonstandard integer.  Why does this produce a character?

\medskip

A.C. - The recipe is very simple to get a character from a nonstandard
integer:
\begin{enumerate}
\item
View an element~$x$ of the compact group~$C^\N$ as a map~$n \to x(n)$
from the integers~$\N$ to the group~$C$ with two elements~$\pm 1$.
\item
Given a non-standard integer~$n$ the evaluation~${}^*x(n)$ gives an
element of~${}^*C$, but since~$C$ is finite one has~${}^*C=C$.
\item
The map~$x \to {}^*x(n)$ is a character of the compact group~$C^\N$
since it is a multiplicative map from~$C^\N$ to~$\pm 1$.
\item
This character cannot be measurable, since otherwise it would be
continuous and hence~$n$ would be standard.
\end{enumerate}

\medskip
 
M.P.S. - I was precisely asking why it is true that, as you mention in
step~(3), the map~$x\to{}^*x(n)$ is a multiplicative map.

\medskip

A.C. - Just because the product~$xy$ of two elements~$x, y$ in the
group~$C^\N$ is defined by the equality~$(xy)(n)=x(n) y(n)$ for all
$n$, and this equality is first order and holds hence also for
non-standard integers.

\medskip

M.P.S. - Then you are using the transfer principle to conclude that we
have an elementary extension?
 
\medskip

A.C. - Yes, I am using the transfer principle%
\footnote{\label{f7}Connes' acknowledgment of his use of the transfer
principle was mentioned in Section~\ref{four} (see
footnote~\ref{tp}).}
to get that if~$z(n)=x(n)y(n)$ for all~$n$ then one has also
${}^*z(n)={}^*x(n)\,{}^*y(n)$ for all non-standard~$n$.

\medskip

M.P.S. - Exploiting the transfer principle presupposes a \emph{model}
where such a principle applies, such as [for example] the ultrapower
one constructed using an ultrafilter, say a selective one.  With such
a model in the background, seeking to exhibit a character in a
canonical fashion would seem to be canonically equivalent to seeking
to exhibit an ultrafilter.  But why not pick the selective one we
started with?%
\footnote{\label{black2}The point about choosing the ultrafilter that
one started with is related to the metaphor of the black box that
canonically returns its input, mentioned in Section~\ref{five} at
footnote~\ref{black1}.}

\medskip

A.C. - 
\end{quote}

Needless to say, \Schutzenberger{} never challenged Connes as above.
However, the exchange is not entirely virtual: it reproduces an
exchange of emails in june 2012, between Connes and the second-named
author.%
\footnote{The email exchange is reproduced here with the consent of
Connes \cite{Co12b}.}
Connes never replied to the last question about ultrafilters (see the
discussion of the forgetful functor at~\eqref{51}).

\section{Definable model of Kanovei and Shelah}
\label{definable}

In 2004, Kanovei and Shelah constructed a definable model of the
hyperreals.  In this section, we explore some of the meta-mathematical
ramifications of their result.

\subsection{What's in a name?}
\label{what}

Let us consider in more detail Connes' comment on \emph{naming} a
hyperreal:
\begin{quote}
What conclusion can one draw about nonstandard analysis?  This means
that, since noone will ever be able \emph{to~name} a nonstandard
number, the theory remains \emph{virtual} (Connes 2001,
\cite[p.~16]{CLS}) [emphasis added--the authors]
\end{quote}
The exact meaning of the verb ``to name'' used by Connes here is not
entirely clear.  Connes provided a hint as to its meaning in 2000, in
the following terms:
\begin{quote}
if you are given a non standard number you can canonically produce a
subset of the interval which is not Lebesgue measurable.  
Now we know
from logic (from results of Paul Cohen and Solovay) that it will
forever be impossible to produce explicitely [sic] a subset of the
real numbers, of the interval~$[0, 1]$, say, that is not Lebesgue
measurable (Connes 2000 \cite[p.~21]{Co99}, \cite[p.~14]{Co00}).
\end{quote}
The hint is the name \emph{Solovay} (Robert M.\ Solovay).  Apparently
Connes is relying on the following result,
which may be found in (Solovay 1970 \cite[p.~3, Theorem 2]{So}).

\begin{theorem}[Solovay (1970, Theorem 2)]
\label{solt}
There is a model\/~$\som$ of set theory ZFC, in which (it is true
that) every set of reals definable from a countable sequence of
ordinals is Lebesgue measurable.
\end{theorem}

\subsection{The Solovay and G\"odel models}
\label{sgm}

The model~$\som$ mentioned in Theorem~\ref{solt} is referred to as
\emph{the Solovay model} by set theorists.  The assumption of
``definability from a countable sequence of ordinals'' includes
definability from a real (and hence such types of definable pointsets
as Borel and projective sets, among others),
since any real can be effectively represented as a countable
sequence of ordinals --- natural numbers, in this case. 

\begin{remark}
\label{42c}
The model~$\som$ contains a submodel~$\somp$ of all sets~$x$ that are
\emph{hereditarily definable from a countable sequence of ordinals}.
This means that~$x$ itself, all elements~$y\in x$, all elements of
elements of~$x$, \emph{etc.}, are definable from a countable sequence
of ordinals.  This submodel~$\somp$ is sometimes called \emph{the
second Solovay model}.  It turns out that~$\somp$ is a model of ZF in
which the full axiom of choice AC fails.  Instead, \emph{the axiom DC
of countable dependent choice}%
\footnote{Given a sequence of nonempty sets~$\langle X_n : n\in \N
\rangle$, the axiom DC postulates the existence of a countable
sequence of choices~$x_0,x_1,x_2,\dots$ in the case when, for
each~$n$, the domain~$X_n$ of the~$n^{\text{th}}$ choice~$x_n\in X_n$
may depend not only on~$n$ but also on the previously made
choices~$x_0,x_1,\dots,x_{n-1}$.  It is considered to be the strongest
possible version of ``countable choice''.}
holds in~$\somp$, so that~$\somp$ is a model of ZF+DC.  
\end{remark}

The following is an immediate consequence of Theorem~\ref{solt}.

\begin{corollary}[Solovay (1970, Theorem 1)]
\label{solt2}
It is true in the second Solovay model\/~$\somp$ that every set of
reals is Lebesgue measurable.
\end{corollary}

\begin{remark}
\label{dc}
A free ultrafilter on~$\N$ yields a set in~$(\Z/2\Z)^\N$ which is
nonmeasurable in the sense of the natural uniform probability measure
on~$(\Z/2\Z)^\N$.  Meanwhile, the second Solovay model~$\somp$ of
ZF+DC contains no such sets, and therefore no such ultrafilters,
either.  It follows that one cannot prove the existence of a free
ultrafilter on~$\N$ in ZF+DC.
\end{remark}

The \emph{constructible model\/}~$\bL$, introduced by (G\"odel 1940,
\cite{go}), is another model of ZFC, opposite to the Solovay model in
many of its features, including the existence of \emph{definable
non-measurable} sets of reals.  Indeed, it is true in~$\bL$ that there
is a non-measurable set in~$\R$ which is not merely definable, but
definable in a rather simple way which places it in the effective
class~$\varDelta^1_2$ of the projective hierarchy (see P.~Novikov
\cite{No}).  With these two models in mind, it is asserted that the
existence of a definable Lebesgue non-measurable set is
\emph{independent} of the axioms of set theory.

\subsection{That which we call a non-sequitur}
\label{what2}

If, in Connes' terminology, ``to name'' is ``to define'', then Connes'
remark to the effect that
\begin{quote}
since noone will ever be able to name a nonstandard \emph{number}, the
\emph{theory} remains virtual (Connes 2001, \cite[p.~16]{CLS})
[emphasis added--the authors]
\end{quote}
is that which we call a \emph{non-sequitur}.  Namely, while an
ultrafilter (associated with a non-standard number by means of the
transfer principle) cannot be \emph{defined}, a \emph{definable}
(countably saturated) model of the hyperreals was constructed by
Kanovei and Shelah (2004, \cite{KS}).  Their construction appeared
later than Connes' ``virtual'' comment cited at the beginning of
Subsection~\ref{what}.  However, three years after the publication of
\cite{KS}, Connes again came back to an alleged ``catch in the
theory'':
\begin{quote}
I had been working on non-standard analysis but after a while I had
found a catch in the theory\ldots The point is that as soon as you
have a non-standard number, you get a non-measurable set.  And in
Choquet's circle, having well studied the Polish
%
%
school, we knew that every set you can \emph{name} is measurable
(Connes 2007, \cite[p.~26]{Co07a}) [emphasis added--the authors].
\end{quote}
An ultrafilter associated with a non-standard \emph{number} cannot be
``named'' or, more precisely, \emph{defined}; however, the
\emph{theory} had been shown (three years prior to Connes' 2007
comment) to admit a \emph{definable} model.  Connes' reference to
Solovay suggests that, to escape being \emph{virtual}, a theory needs
to have a \emph{definable} model.  If so, his ``virtual'' allegation
concerning non-standard analysis is erroneous, by the result of
Kanovei and Shelah.

Connes' claim that ``every set you can name is measurable'' is
similarly inaccurate, by virtue of the G\"odel constructible
model~$\bL$, as discussed in Subsection~\ref{sgm}.  A correct
assertion would be the following: if you ``name'' a set of reals then
you cannot prove (in ZFC) that it is nonmeasurable, and moreover, one
can ``name'' a set of reals (a G\"odel counterexample) regarding which
you cannot prove that it is measurable, either.

Connes elaborated a distinction between countable AC and uncountable
AC, and criticized NSA for relying on the latter (see
Section~\ref{contrast}).  He invoked the Solovay model to explain why
he feels NSA is a ``virtual'' theory.  Now the second Solovay model
${\mathcal S}'$ of ZC+DC demonstrates that ultrafilters on~$\N$ cannot
be shown to exist without uncountable AC (see Remark~\ref{42c}).
Thus, no ultrafilters, chimerical or otherwise, can be produced by
means of the countable axiom of choice alone; yet Connes exploited
ultraproducts (and ultrafilters on~$\N$) in an essential manner in his
work on the classification of factors (Connes 1976, \cite{Co76}).

\section{Machover's critique}
\label{machover}

In 1993, M.~Machover analyzed non-standard analysis and its role in
teaching, expanding on a discussion in J.~Bell \& Machover
\cite[p.~573]{BM}.  We will examine Machover's criticism in this
section.

\subsection{Is there a best enlargement?}
In 1993, Machover wrote:
\begin{quote}
The [integers, rationals, reals] can be characterized (informally or
within set theory) \emph{uniquely up to isomorphism} by virtue of
their mathematical properties \ldots But there is no \ldots known way
of singling out a particular enlargement that can plausibly be
regarded as canonical, nor is there any reason to be sure that a
method for obtaining a canonical enlargement will necessarily be
invented (Machover 1993 \cite{Ma93}) [emphasis added--the authors]
\end{quote}
The problem of the uniqueness of the nonstandard real line is
discussed in detail in an article by Keisler (1994, \cite{keis}), to
which we refer an interested reader.  Meanwhile, Machover emphasizes
\begin{enumerate}
\item[(A)] the \emph{uniqueness up to isomorphism} of the traditional
number systems (integers, rationals, reals), allegedly \emph{unlike}
the hyperreals; and
\item[(B)] an absence of a preferred enlargement.
\end{enumerate}
As we will see, he is off-target on both points (though the latter
became entirely clear only after his text was published).  We start
with three general remarks.  

(1) A methodological misconception on the part of some critics of NSA
is an insufficient appreciation of the fact that the hyperreal
approach does \emph{not} involve a claim to the effect that
hyperreals~$\IR$ are ``better'' than~$\R$.  Rather, one works with the
\emph{pair}~$(\R,\IR)$ together with, say, the standard part function
from limited hyperreals to~$\R$.  It is the interplay of the
\emph{pair} that bestows an advantage on this approach.  The real
field is still present in all its unique complete Archimedean totally
ordered glory.

(2) Noone would dismiss an algebraic number field on the grounds that
it is not as good as~$\Q$ because of a lack of uniqueness.  It goes
without saying that the usefulness of an algebraic number field is not
impaired by the fact that there exist other such number fields.

(3) The specific technical criticism of Machover's that the hyperreal
enlargement is not unique and therefore one needs to prove that the
notion of ``continuity", for example, is model-independent, is
answered by the special enlargement constructed by Morley and Vaught
(1962, \cite{mv}) for any uncountable cardinality~$\kappa$ satisfying
$2^\alpha\le \kappa$ for all~$\alpha<\kappa$ (see Subsection~\ref{AP}
for more details) and providing a unique such enlargement up to
isomorphism.  

\begin{remark}
\label{51b}
Under the assumption of GCH, the condition on~$\kappa$ holds for all
infinite cardinals~$\kappa$.  If GCH is not assumed, then it still
holds for \emph{unboundedly many} uncountable cardinals%
\footnote{Namely, for every cardinal there is one of this kind (note
that this is more than merely ``infinitely many'').}
$\kappa$, one of which (not necessarily the least one) can be defined
by~$\kappa = lim_n a_n$, where~$a_0 = \aleph_0$ and~$a_{n+1} =
2^{a_n}$.
\end{remark}

\subsection{Aesthetic and pragmatic criticisms}
\label{AP}

Machover's critique of NSA actually contains two separate criticisms
even though he tends to conflate the two.  The first criticism is an
\emph{aesthetic} one, mainly addressed to traditionally trained
mathematicians: the reals are unique up to isomorphism, the hyperreals
aren't.  The second criticism is a \emph{pragmatic} one, and is
addressed to workers in NSA: hyperreal definitions of standard
concepts apparently depend on the particular extension of~$\R$ chosen,
and therefore necessitate additional technical work.  We will comment
further on the two criticisms below.

Machover expressed his \emph{aesthetic} criticism by noting that if we
choose a system of real numbers
\begin{quote}
in which the Continuum Hypothesis holds, and another in which it does
not [hold], then for each such choice there are still infinitely many
non-isomorphic enlarged systems of [hyper]reals, none of which has a
claim to be `the best one' (Machover 1993, \cite[p.~210]{Ma93}).
\end{quote}

How cogent is Machover's aesthetic criticism?  The CH-part of his
claim
%
%
is dubious as it does not accord with what we observed above.  Indeed,
as noted in Subsection~\ref{52}, all models of~$\dR$ of the
form~$\R^\N/\mathcal{F}$ are isomorphic in ZFC+CH (see
Theorem~\ref{erdos}).  The uniqueness of the isomorphism type of such
a hyperreal field parallels that of the traditional structures
(integers, rationals, reals) emphasized by Machover in item~(A) above.

The non-CH part of Machover's claim is similarly dubious.  Although
all models of~$\dR$ are not necessarily isomorphic under the ZFC
axioms, still uniqueness up to isomorphism is attainable within the
category of \emph{special} models, that is, those represented in the
form of limits of certain increasing transfinite sequences of
successive saturated elementary extensions of~$\dR$.  (See a detailed
definition in Chang and Keisler (1990, \cite{ck}), 5.1.)  The
following major theorem is due to Morley and Vaught (1962, \cite{mv}),
see also 5.1.8 and 5.1.17 in Chang and Keisler (1990, \cite{ck}).

\begin{theorem}
\label{mvck}
Suppose that an uncountable cardinal\/~$\kappa$ satisfies the
implication\/~$\alpha<\kappa\Longrightarrow 2^\alpha\le\kappa$.  Then
\begin{enumerate}
\item
there are special models of\/~$\dR$ of cardinality\/~$\kappa$, and
\item
all those models are pairwise isomorphic.
\end{enumerate}
\end{theorem}

Thus, for any cardinal~$\kappa$ as in the theorem, there is a uniquely
defined isomorphism type of nonstandard extensions of~$\R$ of
cardinality~$\kappa$.  Cardinals of this type do exist independently
of GCH (see Remark~\ref{51b}) and can be fairly large, but at any rate
one does have uniquely defined isomorphism types of models of~$\IR$
in suitable infinite cardinalities.

\begin{remark}
A decade after the publication of Machover's article, Kanovei and
Shelah (2004, \cite{KS}) proved the existence of a definable
individual model of the hyperreals (not just a definable isomorphism
type), contrary to all expectation (including Machover's, as the
passage cited above suggests).  Further research by Kanovei and
Uspensky (2006, \cite{ku}) proved that all Morley--Vaught isomorphism
classes given by Theorem~\ref{mvck} likewise contain definable
individual models of~$\dR$.%
\footnote{We note that a \emph{maximal} class hyperreal field (in the
von Neumann-Bernays-G\"odel set theory) was recently analyzed by
Kanovei and Reeken (2004, \cite[Theorem 4.1.10(i)]{kr}) in the
framework of axiomatic nonstandard analysis, and by P.~Ehrlich (2012,
\cite{Eh12}) from a different standpoint.  In each version, it is
similarly unique, and, in the second version, isomorphic to a maximal
surreal field.}
\end{remark}

\begin{remark}
If one works in the Solovay model~$\som$ 
as a background ZFC universe, then the definable Kanovei--Shelah model
of~$\dR$ does not contain a \emph{definable nonstandard integer}, as
any such would imply a definable non-measurable set, contrary to
Theorem~\ref{solt}.  The apparent paradox of a non-empty definable set
with no definable element is an ultimate expression of a known
mathematical phenomenon when a simply definable set has no equally
simply definable elements.%
\footnote {For instance, one can define in a few lines what a
transcendental real number is, but it would require a number of pages
to prove for an average math student that~$\pi$,~$e$, or any other
favorite trancendental number is in fact trancendental.}
\end{remark}

As to Machover's \emph{pragmatic} criticism addressed to NSA workers,
we note that requiring suitable properties of saturation in a given
cardinal, one in fact does obtain a unique model of the hyperreals.
Therefore the criticism concerning the dependence on the model becomes
moot.

\subsection{Microcontinuity}
\label{mc}

Machover recalls a property of a function~$f$ that we will refer to as
\emph{microcontinuity} at a point~$r\in\R$ following
Davis~\cite[p.~96]{Da77}:
\begin{equation}
\label{mach}
f(x)\approx f(r) \text{\;for every hyperreal\;} x\approx r.
\end{equation}
Here ``$\approx$'' stands for equality up to an infinitesimal.
Property~\eqref{mach} is equivalent%
\footnote{Strictly speaking~$f$ should be replaced by~${}^*\!f$
in~\eqref{mach}.  Note that, modulo replacing the term ``hyperreal''
by the expression ``variable quantity'', definition~\eqref{mach} is
Cauchy's definition of continuity, contrary to a widespread
Cauchy--Weierstrass tale concerning Cauchy's definition (see Borovik
et al.~\cite{BK} as well as \cite{KK11a, KK11b}).}
to the usual notion of continuity of a real function~$f$ at~$r$.
Machover goes on to assert that
\begin{quote}
in order to legitimize [\eqref{mach}] as a definition \ldots, we must
make sure that it is independent of the choice of
enlargement. (Otherwise, what is being defined would be a
\emph{ternary} relation between~$f$,~$r$ and the enlargement.)
(Machover 1993, \cite[p.~208]{Ma93}).
\end{quote}
Microcontinuity formally depends on the enlargement.  Machover
concludes that it cannot replace~$(\epsilon,\delta)$ definitions
altogether:
\begin{quote}
Therefore, [\eqref{mach}] cannot displace the old standard
[$\epsilon,\delta$] definition altogether, if one's aim is to achieve
proper rigour and methodological correctness \ldots There is a long
tradition of teaching \emph{first-year calculus} in a way that
sacrifices a certain amount of rigour in order to make the material
more intuitive.  There is, of course, nothing wrong or dishonourable
about this--provided the students are told that what they are getting
is a version that does not satisfy the highest standards of rigour and
glosses over some problems requiring closer consideration (ibid.)
[emphasis added--the authors].
\end{quote}
Granted, we need to be truthful toward our students.  However,
Machover's argument is unconvincing, as he misdiagnozes the
educational issue involved.  The issue is \emph{not} whether
the~$(\epsilon,\delta)$ definition should be replaced
\emph{altogether} by a microcontinuous definition as in~\eqref{mach}.
Rather, the issue revolves around \emph{which} definition should be
the \emph{primary} one.  Thus, Keisler's textbook does present
the~$(\epsilon,\delta)$ definition (Keisler 1986, \cite[p.~286]{Ke}),
once continuity has been thoroughly explained via microcontinuity.%
\footnote{Pedagogical advantages of microcontinuity were discussed in
B\l aszczyk et al.~\cite[Appendix A.3]{BKS}.}
The~$(\epsilon,\delta)$ definition is an elementary formula, which
shows that continuity is expressible in first order logic, a fact not
obvious from the microcontinuous definition~\eqref{mach} dependent as
it is on an external relation~``$\approx$''.  Since
the~$(\epsilon,\delta)$ definition needs to be mentioned in any case,
the apparent dependence of~\eqref{mach} on the choice of an
enlargement is a moot point.

\section{How powerful is the transfer principle?}
\label{seven}

The back cover of the 1998 hyperreal textbook by R.~Goldblatt
describes non-standard analysis as
\begin{quote}
a wellspring of powerful new principles of reasoning (transfer,
overflow, saturation, enlargement, hyperfinite approximation, etc.)
(see Goldblatt 1998 \cite{Go}).
\end{quote}
Of the examples mentioned here, we are particularly interested in
\emph{transfer}, i.e., the transfer principle whose roots go back to
\Los's theorem (\Los{} 1955, \cite{Lo55}).  The back cover describes
the transfer principle as a \emph{powerful new principle of
reasoning}.

On the other hand, a well-established tradition started by P.~Halmos
holds that the said principle is not powerful at all.  Thus, Halmos
described non-standard analysis as 
\begin{quote}
a special tool, too special, and other tools can do everything it
does (Halmos 1985, \cite[p.~204]{Ha}).
\end{quote}

Are we to conclude that the 1998 back cover contains a controversial
assertion and/or a well-meaning exaggeration?  Hardly so.  The term
``powerful'' is being used in different senses.  In this section we
will try to clarify some of the meanings of the term.

\subsection{Klein--Fraenkel criterion}
\label{KF}

In 1908, Felix Klein formulated a criterion of what it would take for
a theory of infinitesimals to be successful.  Namely, one must be able
to prove a mean value theorem (MVT) for arbitrary intervals, including
infinitesimal ones:
\begin{quote}
The question naturally arises whether [\ldots] it would be possible to
modify the traditional foundations of infinitesimal calculus, so as to
include actually infinitely small quantities in a way that would
satisfy modern demands as to rigor; in other words, to construct a
non-Archimedean system.  The first and chief problem of this analysis
would be to prove the mean-value theorem
\[
f(x+h)-f(x)=h \cdot f'(x+\vartheta h)
\]
from the assumed axioms.  I will not say that progress in this
direction is impossible, but it is true that none of the investigators
have achieved anything positive (Klein 1908, \cite[p.~219]{Kl08}).
\end{quote}
In 1928, A.~Fraenkel \cite[pp.~116-117]{Fran} formulated a similar
criterion in terms of the MVT.

Such a Klein--Fraenkel criterion is satisfied by the framework
developed by Hewitt, \Los, and Robinson.  Indeed, the MVT is true for
the natural extension~${}^*\!f$ of every real smooth function~$f$ on
an arbitrary hyperreal interval, by the transfer principle.
Fraenkel's opinion of Robinson's theory is on record:
\begin{quote}
my former student Abraham Robinson had succeeded in saving the honour
of infinitesimals - although in quite a different way than Cohen%
\footnote{The reference is to Hermann Cohen (1842--1918), whose
fascination with infinitesimals elicited fierce criticism by both
G.~Cantor and B.~Russell.  For an analysis of Russell's
\emph{non-sequiturs}, see Ehrlich \cite{Eh06} and Katz \& Sherry
\cite{KS1}, \cite{KS2}.}
and his school had imagined (Fraenkel 1967, \cite[p.~107]{Fra67}).
\end{quote}

The hyperreal framework is the only modern theory of infinitesimals
that satisfies the Klein-Fraenkel criterion.  The fact that it
satisfies the criterion is due to the transfer principle.  In this
sense, the transfer principle can be said to be a ``powerful new
principle of reasoning''.

One could object that the classical form of the MVT is not a key
result in modern analysis.  Thus, in L. H\"ormander's theory of
partial differential operators (H\"ormander \cite[p.~12--13]{Ho}), a
key role is played by various multivariate generalisations of the
following Taylor (integral) remainder formula:
\begin{equation}
\label{71}
f(b)=f(a)+(b-a)f'(a)+\int_a^b (b-x)f''(x) dx.
\end{equation}
Denoting by~$\mathcal{D}$ the differentiation operator and by
$\mathcal{I}=\mathcal{I}(f,a,b)$ the definite integration operator, we
can state~\eqref{71} in the following more detailed form for a
function~$f$:
\begin{equation}
\label{72}
\begin{aligned}
(\forall & a\in\R) (\forall b\in\R) \\ & f(b)
=f(a)+(b-a)(\mathcal{D}f)(a)+ \mathcal{I} \left(
(b-x)(\mathcal{D}^2_{\phantom{I}} f),a,b \right)
\end{aligned}
\end{equation}
Applying the transfer principle to the elementary formula~\eqref{72},
we obtain
\begin{equation}
\label{73}
\begin{aligned}
(\forall & a\in\IR) (\forall b\in\IR) \\ &
{}^*\!f(b)={}^*\!f(a)+(b-a)({}^*\mathcal{D}\;{}^*\!f)(a)+
{}^*\mathcal{I} \left( (b-x)({}^*\mathcal{D}^2_{\phantom{I}}\;
{}^*\!f),a,b \right)
\end{aligned}
\end{equation}
for the natural hyperreal extension~${}^*\!f$ of~$f$.  The
formula~\eqref{73} is valid on every hyperreal interval of~$\IR$.
Multivariate generalisations of~\eqref{71} can be handled similarly.

We focused on the MVT (and its generalisations) because, historically
speaking, it was emphasized by Klein and Fraenkel.  The transfer
principle applies far more broadly, as can be readily guessed from the
above.

\subsection{Logic and physics}
\label{agree}

There is another sense of the term \emph{powerful} that is more
controversial than the one discussed in Subsection~\ref{KF}.  Namely,
how powerful are the hyperreals as a research tool and an engine of
discovery of \emph{new} mathematics?  The usual litany of impressive
breakthroughs achieved using NSA includes progress on the invariant
subspace problem, canards, hydrodynamics and Boltzmann equation,
non-standard proof of Gromov's theorem on groups of polynomial growth,
Hilbert's fifth problem (see Hirschfeld \cite{Hi} and Goldbring
\cite{Go10}), etc.%
\footnote{For additional examples see the book \cite{VN}.}

However, declaiming such a list does little more than encourage the
partisans while further antagonizing the critics.  We will therefore
comment no further other than clarifying that this is \emph{not} the
meaning of the term \emph{powerful} when we use it in reference to the
transfer principle.  Namely, we use it solely in the sense explained
in Subsection~\ref{KF}.

The significance of the back cover comment cited at the beginning of
Section~\ref{seven} is that Robinson's theory introduces new
perspectives and intuitions into mathematics, similarly to physics.%
\footnote{Such an analogy between logic and physics is due to David
Kazhdan.}
When E.~Witten informally wrote down a pair of equations on the board
at MIT a couple of decades ago, he was motivated by physical
intuitions.  The resulting Seiberg-Witten theory caused a revolution
in gauge theory, and in particular resulted in much shorter proofs of
theorems that S.~Donaldson received his Fields medal for (see e.g.,
Katz \cite{Ka}).  Logic, similarly, introduces new intuitions and
techniques.  Today logicians like E.~Hrushovski \cite{Hru} obtain
results in ``ordinary mathematics" by model-theoretic means.

Interesting recent uses of non-standard methods as applied to the
structure of approximate groups may be found in Hrushovski
\cite{Hru12} and Breuillard, Green, \& Tao \cite{BGT12}.

\section{How non-constructive is the Dixmier trace?}
\label{dix}

This section deals with the foundational status of the Dixmier trace,
and with the role of Dixmier trace in noncommutative geometry.

\subsection{Front cover}
\label{cover}

The front cover of the book \emph{Noncommutative geometry} features an
elaborate drawing, done by Connes himself (according to the copyright
page).  The drawing contains only three formulas.  One of them is the
expression
\[
-\hskip-13pt\int |dZ|^p.
\]
The barred integral symbol~$-\hskip-11pt\int$ is Connes' notation for
the trace constructed by Dixmier (1966 \cite{Di66}).  The notation
first occurred in print in (Connes 1995 \cite[p.~6213, formula
(2.34)]{Co95}), i.e., the year after its appearance on the front cover
of Connes' book.  The appearance of Dixmier's trace on the book cover
indicates not only that Connes was already thinking of the Dixmier
trace as a kind of ``integration'' (this idea is already found in
Connes 1988 \cite{Co88}), but also that Connes himself thought of the
trace as an important ingredient of noncommutative geometry.

\subsection{Foundational status of Dixmier trace}

The Dixmier trace is a linear functional on the space of compact
operators whose characteristic values have a specific rate of
convergence to~$0$.  In Connes' framework, the Dixmier trace can be
thought of as a kind of an ``integral'' of infinitesimals.  An
analogous concept in Robinson's framework is the functional
\[
\text{st}(n\epsilon).
\]
Here~$n\in{}^*\N\setminus\N$ is a \emph{fixed} hypernatural, and the
functional is defined for a \emph{variable} infinitesimal~$\epsilon$
constrained by the condition that~$n\epsilon$ is finite.

Dixmier exploited ultrafilters in constructing his trace.  Dixmier
traces can also be constructed using universally measurable
\emph{medial limits}, independently constructed by Christensen
\cite{Chr} and Mokobodzki in the assumption of the continuum
hypothesis (CH).  Mokobodzki's work was explained by P.~Meyer (1973,
\cite{Me}).  Meyer's text is cited in Connes' book \cite{Co94}, but
not in the section dealing with Dixmier traces (Connes 1994,
\cite[p.~303-308]{Co94}), which does not use medial limits and instead
relies on the Hahn--Banach theorem \cite[p.~305, line 8 from
bottom]{Co94}.%
\footnote{Note that, in the spirit of reverse mathematics, the
Hahn--Banach theorem is sufficient to generate a Lebesgue
nonmeasurable set (see Foreman \& Wehrung \cite{FW},
Pawlikowski~\cite{Pa}).}

Medial limits have been shown \emph{not} to exist in the assumption of
the filter dichotomy (FD) by P.~Larson \cite{La09}.  FD is known to be
consistent (Blass and Laflamme \cite{BL}).  The assumption of CH
(exploited in the construction of medial limits) is generally
considered to be a very strong foundational assumption, more
controversial than the axiom of choice (see e.g., J.~Hamkins
\cite{Ham}, \cite{Ham12}; D.~Isaacson \cite{Is}).

Indeed, while all the major applications of the ``uncountable''\ AC\
outside of set theory proper%
\footnote{This includes such constructions as the Vitali
non-measurable set, Hausdorff's gap, ultrafilters on~$\N$, the Hamel
basis, the Banach--Tarski paradox, nonstandard models, etc.
Sierpi\'nski (1934, \cite{sie}) gives many additional examples.}
can be reduced to the assumption that the continuum of real numbers
can be wellordered, CH requires, in addition, the existence of a
wellordering of~$\R$ specifically of length~$\omega_1$ (which is the
shortest possible length of such a wellordering).

Moreover, CH implies the existence of~$P$-point ultrafilters%
\footnote{See footnote~\ref{rudin}.}
on~$\N$, and Shelah \cite{Sh} showed that the existence of~$P$-points
cannot be established in ZFC, again indicating the controversial
nature of CH.

Furthermore, Connes notes that the results he is interested in happen
to be \emph{independent of the choice} of the Dixmier trace
\cite[p.~307, line 14 from bottom]{Co94}.  Thus the strong assumption
of CH appears superfluous, and the nonconstructive nature of the
ultrafilter construction of the Dixmier trace, a paper tiger.  Namely,
Dixmier trace is constructive or non-constructive in a sense similar
to that of a hyperreal number being constructive or non-constructive:
both rely on nonconstructive foundational material (be it AC, CH, or
Hahn-Banach), but yield results \emph{independent of choices} made.
For instance, differentiating~$x^2$ yields~$2x$ regardless of the
variety of infinitesimals exploited in defining the derivative.
Similarly, the notion of continuity, when defined via microcontinuity,
is independent of the hyperreal model used (see Subsection~\ref{mc}).

\subsection{Role of Dixmier trace in noncommutative geometry}

At a recent conference (see \cite{GIS}) on singular traces (such as
the Dixmier trace), a majority of the speakers mentioned the Dixmier
trace in their abstracts, while none of them mentioned (or cited)
either Mokobodzki or medial limits.  Recent work by the conference
speakers dealing with Dixmier traces includes: Carey, Phillips, \&
Sukochev~\cite{CPS}; Engli\v s \& Zhang \cite{EZ}; Lord \& Sukochev
\cite{LS10, LS11}; Lord, Potapov, \& Sukochev \cite{LPS}; Kalton,
Sedaev, \& Sukochev \cite{KSS}; Sukochev \& Zanin \cite{SZ1, SZ2}.

Most speakers also cite Connes' \emph{Noncommutative Geometry}.  Ever
since its appearance on the front cover of Connes' book (see
Subsection~\ref{cover}), the Dixmier trace has played a major role in
Connes' framework and related fields.

\section{Of darts, infinitesimals, and chimeras}
\label{dart}

In this section we will be concerned with a somewhat elusive issue of
what is real and what is chimerical.

\subsection{Darts}
\label{darts}

Connes outlined a game of darts in 2000 in the following terms:

\begin{quote}
You play a game of throwing darts at some target called~$\Omega$
\ldots{} what is the probability~$dp(x)$ that actually when you send
the dart you land exactly at a given point~$x\in\Omega$~? \ldots{}
what you find out is that~$dp(x)$ is smaller than any positive real
number~$\varepsilon$.  On the other hand, if you give the answer that
$dp(x)$ is~$0$, this is not really satisfactory, because whenever you
send the dart it will land \emph{somewhere} (Connes 2000,
\cite[p.~13]{Co00}) [emphasis added--the authors].
\end{quote}
As Connes points out, no satisfactory interpretation of such
intuitions seems to exist in a real number system devoid of
infinitesimals.  But if one interprets the ``$p$" to be an
infinitesimal interval rather than a point, there is a consistent
theory that can capture the intuitions Connes spoke of.  Namely,
assume for the sake of simplicity that the target is the unit
interval~$[0,1]$.  A more satisfactory answer than the one above is
provided in terms of a hyperfinite grid
\begin{equation}
\label{91}
\text{Grid}_H= \left\{ 0, \tfrac{1}{H}, \tfrac{2}{H}, \tfrac{3}{H},
\ldots, \tfrac{H-1}{H}, 1 \right\}
\end{equation}
defined by a hypernatural~$H\in{}^*\N\setminus\N$.  Then the
probability of the dart hitting an infinitesimal
interval~$[\tfrac{k}{H}, \tfrac{k+1}{H}] \subset [0,1]$ can be taken
to be precisely~$\tfrac{1}{H}$.  The hypernatural~$H$ can be chosen to
be the explicit unchimerical one appearing in~\eqref{31}.

Similarly, the probability of the dart hitting a real set
$A\subset[0,1]$ can be computed as follows.  Roughly speaking, one
counts the number of points in the intersection
$\text{st}^{-1}(A)\cap\text{Grid}_H$ and divides by~$H$, where
$\text{st}$ is the standard part function on limited hyperreals, and
$\text{Grid}_H$ is the hyperfinite grid of~\eqref{91}, yielding a
probability of
\begin{equation}
\label{92}
\frac{|\text{st}^{-1}(A)\cap\text{Grid}_H|}{H};
\end{equation}
more precisely, since~$\text{st}^{-1}(A)$ is not an internal set, one
takes the infimum of~$\text{st}(|X|/H)$ over all internal sets~$X$
containing~$\text{st}^{-1}(A)\cap\text{Grid}_H$ (see
Goldblatt~\cite{Go}, Lemma~16.5.1 on page~210, and Theorem~16.8.2 on
page~217).
%
%

\subsection{Chimeras}
\label{chimera}

Probability theory and measure theory over the hyperreals are today
vast research fields (see e.g., Benci et al.~\cite{BHW}, Wenmackers \&
Horsten~\cite{WH}).  Meanwhile, Connes comments as follows:
\begin{quote}
A nonstandard number is some sort of chimera%
\footnote{\label{uccello2}See Figure~\ref{uccello}.  The reader may be
amused to find similar terminology in Karl Marx, who commented as
follows: ``The closely held belief of some rationalising
mathematicians that~$dy$ and~$dx$ are quantitatively actually only
infinitely small, only approaching 0/0, is a chimera'' (Marx cited in
Fahey et al.~\cite[p.~260]{Fa09}).}
which is impossible to grasp and certainly not a concrete object.  In
fact when you look at nonstandard analysis you find out that except
for the use of ultraproducts, which is \emph{very efficient},
it just shifts the order in logic by one step; it's not doing much
more (Connes 2000, \cite[p.~14]{Co00}) [emphasis added--the authors]
\end{quote}

Connes describes ultraproducts as ``very efficient'', apparently in
contrast to the rest of non-standard analysis.  Meanwhile, the special
case of an ultraproduct used in the construction of~$\IR$ as
in~\eqref{42} exploits an ultrafilter~$\mathcal{F}$ described by
Connes as a ``chimera''.  Are we to conclude that we are dealing with
a \emph{very efficient chimera}?

\begin{remark}
\label{von}
Connes exploits a nonprincipal ultrafilter~$\omega$ in constructing
the ultraproduct von Neumann algebra~$N^\omega$ containing a von
Neumann algebra~$N$ in \emph{Noncommutative geometry}:
\begin{quote}
Definition 11. For every ultrafilter~$\omega\in\beta\N\setminus\N$ let
$N^\omega$ be the ultraproduct,~$N^\omega=$ the von Neumann algebra
$\ell^\infty(\N,N)$ divided by the ideal of sequences~$(x_n)_{n\in\N}$
such that~$\lim_{n\to\omega} \|x_n\|_2 = 0$ (Connes \cite[ch.~V,
sect.~6.$\delta$, Def.~11]{Co94}).%
\footnote{The definition appears on page 495 in the pdf version
available from Connes' homepage, and on page 483 in the published
book.}
\end{quote}
\end{remark}

Perhaps Connes' intention is similar to that of Leibniz, who sometimes
described infinitesimals as ``useful fictions'' (see Katz \& Sherry
\cite{KS1, KS2} and Section~\ref{otte} below).  But Leibniz's position
is generally thought to be close to a formalist one, akin to
Robinson's, whereas Connes is known as a Platonist (see
Subsection~\ref{103}).%
\footnote{\label{plato1}See also footnote~\ref{plato2} on a comment by
Davies.}

Connes goes on to argue that \emph{his} infinitesimal framework does
provide an adequate framework for solving the dart problem (see
\cite[formula (2.35)]{Co97}).  However, Connes' noncommutative
infinitesimals do not form a division ring, do not possess a total
order, lack a transfer principle, and would have difficulty handling
the dart problem as smoothly as~\eqref{92}.

\subsection{Shift}
\label{shift}

What is the meaning of the phrase
\begin{quote}
``nonstandard analysis \ldots just shifts the order in logic by one
step; it's not doing much more''
\end{quote}
penned by Connes (see Subsection~\ref{chimera})?  The phrase is
characteristically evasive (cf.~the discussion of his use of the verb
``to name'' in Subsection~\ref{what}), but perhaps he is referring to
the fact that non-standard analysis permits one to express formulas in
\emph{second} order logic as formulas in \emph{first} order logic over
the hyperreals (hence ``shifts the order in logic by one step'').  In
this context, it is instructive to consider what Fields medalist
T. Tao has to say concerning the expressive power of non-standard
analysis:
\begin{quote}
[it] allows one to rigorously manipulate things such as ``the set of
all small numbers'', or to rigorously say things like ``$\eta_1$ is
smaller than anything that involves~$\eta_0$'', while greatly reducing
epsilon management issues by automatically concealing many of the
quantifiers in one's argument (Tao 2008 \cite[p.~55]{Tao08}).
\end{quote}
The 2009 Abel prize winner M.~Gromov said in 2010:
\begin{quote}
After proving the theorem about polynomial growth using the limit and
looking from infinity, there was a paper by Van den Dries and Wilkie
giving a much better presentation of this using ultrafilters (Gromov
cited in \cite{RS}).
\end{quote}

Other authors have taken note of Connes' sweeping judgments of
mathematical subjects not to his liking.  Thus, E.~B.~Davies writes:
\begin{quote}
In 2001 Alain Connes, a committed Platonist,%
\footnote{\label{plato2}In the context of Davies' comment on Connes'
Platonism, see also main text at footnote~\ref{plato1} which examines
the possibility that Connes may also hold views close to Formalism.}
who has spent a lifetime working on C*-algebras and their
applications, nevertheless excluded the theory of \emph{Jordan
algebras} from the Platonic world of mathematics \ldots How do
mathematicians make such value judgments, and are their opinions more
than prejudices? (Davies 2011, \cite[p.~1456]{Da11}) [emphasis
added--the authors].
\end{quote}
Here Davies is referring to the following comment by Connes:
\begin{quote}
I would say that the exceptional algebra of three-by-three matrices on
Cayley octonions definitely exists because of its connections to the
Lie group~$\text{F}_4$.  As for the general notion of Jordan algebra,
it is difficult to assert that it really holds water
\cite[p.~30]{CLS}.
\end{quote}
Connes finds it ``difficult to assert'' that the theory of Jordan
algebras ``holds water''.  Meanwhile, E.~Zelmanov wrote that
I.~Kantor's work on Jordan algebras (see, e.g., the influential text
Kantor \cite{Kant})
\begin{quote}
played a crucial role in [Zelmanov's] proof of the Restricted Burnside
problem \cite[p.~111]{Ze08},
\end{quote}
work for which Zelmanov was awarded the Fields medal in 1994.

\subsection{Continuum in quantum theory}
\label{real}

Quantum physicists \v Caslav Brukner and 2010 Wolf prize winner
Anton Zeilinger speculate that
\begin{quote}
the concept of an infinite number of complementary observables and
therefore, indirectly, the assumption of continuous variables, are
just mathematical constructions which might not have a place in a
final formulation of quantum mechanics \ldots continuous variables are
devoid of operational and therefore physical meaning in quantum
mechanics" (Brukner \& Zeilinger \cite[p.~59]{BZ}).
\end{quote}
I.~Durham concurs:
\begin{quote}
This latter proposal%
\footnote{I.e., a proposal to resolve the paradox of quantum behavior
of light.}
is similar to coarse-graining arguments in thermodynamic and quantum
systems which have been used by Brukner and Zeilinger to argue that
the continuum is nothing but a mathematical construct, a view I
wholeheartedly endorse (Durham \cite{Du}).
\end{quote}
In 1994, Wolf prize winner John A. Wheeler wrote:
\begin{quote}
The space continuum?  Even continuum existence itself?  Except as
\emph{idealization} neither the one entity nor the other can make any
claim to be a primordial category in the description of nature
(Wheeler \cite[p.~308]{Whe}) [emphasis added--the authors].
\end{quote}
There appears to be an identifiable view in the quantum physics
community that the mathematical continuum is an idealisation, or to
borrow Connes' terminology, it is a ``virtual theory'' or ``chimera'',
though undoubtedly an ``efficient'' one.

A mathematician need not ordinarily be concerned about opinions found
in a separate scientific community.  However, Connes' motivation for
his framework is drawn from quantum theory, and he frequently mentions
quantum mechanics as the inspiration for his noncommutative solution
of the dart problem (see Subsection~\ref{darts}).  His references to
alleged ``absolutely major flaw'' and ``irredemiable defect'' in
Robinson's infinitesimals emanate from their status as an
idealisation.  But in quantum theory, the same observation would apply
to Connes' framework based as it is on the continuum, creating
tensions with Connes' Platonism about the latter (see
Section~\ref{otte}).

Connes claimed that ``it seemed utterly doomed to failure to try to
use non-standard analysis to do physics'' (Connes 2007,
\cite[p.~26]{Co07a}).  Such a claim is particularly dubious coming as
it does two decades after the publication of the 500-page monograph
\emph{Nonstandard Methods in Stochastic Analysis and Mathematical
Physics} by the 1992 Max-Planck-Award recipient S.~Albeverio and
others~\cite{Al}, where just such applications were developed in great
detail.

\section{Conclusion}

The use of non-constructive foundational material such as the axiom of
choice in the hyperreal context is similar to the use of
non-constructive foundational material in Connes' theory.  Thus,
Connes exploits the Dixmier trace (Connes 1995, \cite[p.~6208]{Co95}),
the Hahn--Banach theorem (Connes 1994, \cite[p.~305]{Co94}), as well
as ultrafilters (Connes 1994, \cite[p.~483]{Co94}, see our
Remark~\ref{von} above).  Such concepts rely on non-constructive
foundational material and are unavailable in the framework of the
Zermelo--Fraenkel axioms alone.

Connes claims to provide ``substantial and calculable'' results based
on his theory exploiting the Dixmier trace \cite[p.~211]{Co97}, and
laments the allegedly non-exhibitable nature of Robinson's
infinitesimals.  Meanwhile, Dixmier's construction of the trace relies
on the choice of a nonprincipal ultrafilter on the integers
\cite{Di66}, while an alternative construction requires the continuum
hypothesis (see Section~\ref{dix}).  Connes exploits ultrafilters in
classifying factors and in constructing von Neumann algebras, but
there are no ultrafilters in the second Solovay model~$\mathcal{S}'$
of the set-theoretic universe ZFC+DC (countable choice only) that
Connes professes to favor.  Connes proclaims himself to be an adherent
of countable AC (see Section~\ref{contrast} above), but~$\mathcal{S}'$
is a model of ZFC+DC containing no ultrafilters, so that Connes'
philosophical advocacy of countable AC is divorced from the facts on
the ground of his scientific practice.

Thus, Connes' claims to the effect that his theory produces
computationally meaningful results, allegedly \emph{unlike} Robinson's
theory, are unconvincing.  There is in fact strong similarity between
the two nonconstructivities involved.

Given powerful%
\footnote{See Section~\ref{seven} for a discussion of the term.}
tools such as non-standard enlargements and the transfer principle,
one is able to associate an ultrafilter to a hyperinteger.  But such
ability is a spin-off of the power of the new principles of reasoning
developed in Robinson's approach, and is a reflection, not of a
shortcoming, but rather of the strength of Robinson's method.

\section*{Acknowledgments}

V. Kanovei is grateful to the Fields Institute for its support during
a visit in 2012.  The research of Mikhail Katz was partially funded by
the Israel Science Foundation grant 1517/12.  The research of Thomas
Mormann for this work is part of the research project FFI 2009--12882
funded by the Spanish Ministry of Science and Innovation.

We are grateful to the referees for numerous insightful suggestions
that helped improve an earlier version of the article.  We thank Piotr
B\l aszczyk, Brian Davies, Martin Davis, Ilijas Farah, Jens Erik
Fenstad, Ian Hacking, Reuben Hersh, Yoram Hirshfeld, Karel
Hrb\'a\v{c}ek, Jerome Keisler, Semen Kutateladze, Jean-Pierre Marquis,
Colin McLarty, Elemer Rosinger, David Sherry, Javier Thayer, Alasdair
Urqu\-hart, Lou van den Dries, and Pavol Zlato\v s for helpful
historical and mathematical comments.  The influence of Hilton Kramer
(1928-2012) is obvious.

\bigskip\noindent \textbf{Vladimir Kanovei} graduated in 1973 from
Moscow State University, and obtained a PhD in physics and mathematics
from Moscow State University in 1976.  In 1986, he became Doctor of
science in physics and mathematics at Moscow Steklov Mathematical
Institute (MIAN). He is currently Leading Researcher at the Institute
for Information Transmission Problems (IPPI), Moscow, Russia.  Among
his publications is the book Borel equivalence relations.  Structure
and classification.  \emph{University Lecture Series} \textbf{44}.
American Mathematical Society, Providence, RI, 2008. x+240 pp.

\medskip\noindent \textbf{Mikhail G. Katz} (B.A. Harvard University,
'80; Ph.D. Columbia University, '84) is Professor of Mathematics at
Bar Ilan University.  Among his publications is the book Systolic
geometry and topology, with an appendix by Jake P. Solomon,
\emph{Mathematical Surveys and Monographs}, \textbf{137}, American
Mathematical Society, Providence, RI, '07; and the article (with
A.~Borovik and R.~Jin) An integer construction of infinitesimals:
Toward a theory of Eudoxus hyperreals, \emph{Notre Dame Journal of
Formal Logic} \textbf{53} ('12), no.~4, 557-570.

\medskip\noindent \textbf{Thomas Mormann} studied mathematics,
linguistics and philosophy and earned his PhD in mathematics at the
University of Dortmund in 1978.  Later, he obtained his habilitation
in philosophy at the University of Munich (Germany).  Currently he is
professor of philosophy at the University of the Basque Country
(UPV/EHU) in Donostia-San Sebasti\'an (Spain).  His main fields of
research are philosophy of science and formal ontology.  Among his
publications are: A Place for Pragmatism in the Dynamics of Reason?,
\emph{Studies in History and Philosophy of Science} 43(1), 27--37,
2012; and On the Mereological Structure of Complex States of Affairs,
\emph{Synthese} 187(2), 403--418, 2012.  

\end{document}